\documentclass[11pt]{article}
\usepackage{pxfonts}
\usepackage{yfonts}
\usepackage{dsfont}
\usepackage{graphicx}
\usepackage{relsize}

\def\bel{\begin{equation}\label}
\def\eeq{\end{equation}}
\def\ds{\displaystyle}

\def\mt{\longrightarrow}
\def\v{\vskip 1em}

\def\R{\mathbb R}
\def\Z{\mathbb Z}
\def\C{\mathfrak{C}}

\def\S{{\bf S}}

\def\F{\mathfrak{F}}

\def\O{{\bf O}}
\def\Q{\mathfrak{Q}}

\def\J{{\bf J}}

\def\B{{\bf B}}
\def\H{{\bf H}}
\def\L{{\bf L}}

\def\BMO{{\bf BMO}}

\def\p{{\partial}}

\def\i{{\bf i}}

\def\Hat{\widehat}

\def\bar{\overline}
\def\supp{{\bf supp}}

\def\I{{\bf I}}

\def\M{{\bf M}}

\def\Cup{{\bigcup}}
\def\Cap{{\bigcap}}

\def\alpha{\alphaup}
\def\beta{\betaup}
\def\gamma{\gammaup}
\def\delta{\deltaup}
\def\xi{{\xiup}}
\def\eta{{\etaup}}
\def\tau{{\tauup}}
\def\rho{{\rhoup}}
\def\phi{{\phiup}}
\def\psi{{\psiup}}
\def\lambda{{\lambdaup}}
\def\omega{\omegaup}
\def\varphi{{\varphiup}}
\def\gamma{{\gammaup}}
\def\c{{\bf c}}

\newtheorem{remark}{Remark}[section]

\begin{document}
 \[\begin{array}{cc}\hbox{\LARGE{\bf Regularity of multi-parameter Fourier integral operator}}
 \end{array}\]

 \[\hbox{Zipeng Wang}\]

 \begin{abstract}
We study a family of Fourier integral operators by allowing their symbols to satisfy  a multi-parameter  differential inequality. We  extend the sharp $\L^p$-theorem obtained by Seeger, Sogge and Stein to product spaces. 
\end{abstract}
\section{Introduction}
\setcounter{equation}{0}
In this paper,  we consider a Fourier integral operator defined by
\bel{Ff}
\begin{array}{cc}\ds
\F f(x)~=~\int_{\R^n} f(y)\Omega(x,y)dy,
\\\\ \ds
\Omega(x,y)~=~\int_{\R^n} e^{2\pi\i (\Phi(x,\xi)-y\cdot\xi )}\sigma(x,\xi)d\xi.
\end{array}
\eeq
The symbol  function 
$\sigma(x,\xi)\in\mathcal{C}^\infty(\R^n\times\R^n)$ has a compact support in  $x$. 
On the other hand, the phase function $\Phi(x,\xi)$ is  real, homogeneous of degree $1$ in $\xi$ and smooth for every $x$. Moreover, it satisfies   
 the non-degeneracy condition
\bel{nondegeneracy}
 \det\left[{\p^2\Phi\over \p x\p \xi}\right]\left(x,\xi\right)~\neq~0
 \eeq
at $\xi\neq0$ on the support of $\sigma(x,\xi)$.  

For more background of $\F$, we refer to  the classical reference by Sogge \cite{Sogge}.

$\diamond$  {\small Throughout, we regard $\C$ as a generic constant depending on its subindices.}

We say $\sigma\in S^m$ if 
\bel{class}
\left|\p_\xi^\alpha\p_{x}^\beta \sigma(x,\xi)\right|~\leq~\C_{\alpha~\beta}\left(1+|\xi|\right)^m \left({1\over 1+|\xi|}\right)^{|\alpha|}
\eeq
for every multi-indices $\alpha,\beta$.
 
For $\sigma\in S^0$, 
$\F$ defined as (\ref{Ff})-(\ref{nondegeneracy}) is  bounded on $\L^2(\R^n)$ as shown by  Eskin \cite{Eskin} and H\"{o}rmander \cite{Hormander}. In contrast to this $\L^2$-result,  it is well known that $\F$ of order zero is not bounded on $\L^p(\R^n)$ if $p\neq2$. 

The optimal $\L^p$-estimate was first investigated by   Duistermaat and H\"{o}rmander  \cite{Duistermaat-Hormander} and  then by  Colin de Verdi\'{e}re and Frisch \cite{Colin-Frisch},  Brenner \cite{Brenner}, Peral \cite{Peral}, Miyachi \cite{Miyachi}, Beals \cite{Beals} and eventually obtained by Seeger, Sogge and Stein \cite{S.S.S}.

 \v

{\bf Theorem One: ~Seeger, Sogge and Stein, 1991}\\
 {\it Let $\F$ defined as (\ref{Ff})-(\ref{nondegeneracy}). Suppose $\sigma\in \hbox{S}^m$ for $-(n-1)/2<m\leq0$. We have
\[\left\| \F f\right\|_{\L^p(\R^n)}~\leq~\C_{p~\sigma~\Phi}~\left\| f\right\|_{\L^p(\R^n)},\qquad 1<p<\infty
\]
whenever
\[ \left| {1\over 2}-{1\over p}\right|~\leq~{-m\over n-1}.\]}

{\bf Remark One}~~{\it
This result is sharp. Consider $a(x)b(y)\in\mathcal{C}^\infty_o(\R^n\times\R^n)$ where $a(x)\neq0$ for $|x|=1$ and $b(y)\equiv1$ for $|y|<1$.  Define
\bel{example}
\sigma(x,y,\xi)~=~ a(x)b(y)\left(1+|\xi|\right)^m,\qquad \Phi(x,\xi)~=~ x\cdot\xi+|\xi|.
\eeq
 Then $\F$ given by (\ref{Ff}) is not bounded on $\L^p(\R^n)$ if $\left|1/2-1/p\right|>-m/(n-1),~(1-n)/2\leq m\leq0$.  Regarding estimates can be found at {\bf 6.13}, chapter IX in the book of Stein \cite{Stein}.}

Now,  define $\sigma\in\S^m$ if
\bel{Class}
\left|\p_\xi^\alpha\p_{x}^\beta \sigma(x,\xi)\right|~\leq~\C_{\alpha~\beta}~\left(1+|\xi|\right)^m\prod_{i=1}^n \left({1\over 1+|\xi_i|}\right)^{\alpha_i}
\eeq
for every multi-indices $\alpha,\beta$. 

We give an extension of {\bf Theorem One} by considering the Fourier integral operator $\F$ with a symbol $\sigma\in\S^m$ satisfying the differential inequality in (\ref{Class}).
The study of such operators  that  commute with a multi-parameter family of dilations  dates back to the time of  Jessen, Marcinkiewicz and Zygmund.  Over the several past  decades,
 a number of pioneering  results  have been accomplished, for example   
by Robert Fefferman \cite{R.Fefferman}, Fefferman and Stein \cite{R-F.S}, Chang and  Fefferman \cite{Chang-Fefferman}, Cordoba and Fefferman \cite{Cordoba-Fefferman} 
and M\"{u}ller, Ricci and Stein \cite{M.R.S}. 
\v

{\bf Theorem Two}~~
 {\it Let $\F$ defined as (\ref{Ff})-(\ref{nondegeneracy}). Suppose $\sigma\in \S^m$ for $-(n-1)/2<m\leq0$. We have
  \bel{RESULT}\left\| \F f\right\|_{\L^p(\R^n)}~\leq~\C_{p~\sigma~\Phi}~\left\| f\right\|_{\L^p(\R^n)},\qquad 1<p<\infty
 \eeq
 whenever
 \bel{FORMULA}
 \left| {1\over 2}-{1\over p}\right|~\leq~{-m\over n-1}.
\eeq}

In the next section, we sketch the proof of {\bf Theorem Two} by developing a new framework where the frequency space is decomposed into an infinitely many dyadic cones. 
Every consisting partial  operator whose    symbol  has a support in one of these dyadic  cones is  essentially   an  one-parameter Fourier integral operator, satisfying the desired regularity.

\section{Cone decomposition on frequency space}
\setcounter{equation}{0}
We write $\xi=(\tau,\lambda)\in\R\times\R^{n-1}$. By symmetry, it is suffice to consider  $|\tau|=\max_{i=1,2,\ldots,n}|\xi_i|$. 
Without loss of generality, fix $\tau=\xi_n$. Moreover, 
$\ell$ denotes an $(n-1)$-tuple $(\ell_1,\ell_2,\ldots,\ell_{n-1})$ where  $\ell_i, i=1,2,\ldots,n-1$ are non-negative integers. 

Let $\varphi$ be a smooth {\it bump}-function on $\R$ such that 
\bel{varphi}
\varphi(t)~=~1~~~~\hbox{for}~~~~ |t|~\leq~1,\qquad \varphi(t)~=~0~~~~\hbox{for} ~~~~|t|~>~2.
\eeq
 Define
\bel{delta_t}
\begin{array}{cc}\ds
 \delta_\ell(\xi)~=~\prod_{i=1}^{n-1}\delta_{\ell_i}(\xi),
 \\\\ \ds
 \delta_{\ell_i}(\xi)~=~\varphi\left(2^{\ell_i}{\lambda_i\over\tau}\right)-\varphi\left(2^{\ell_i+1}{\lambda_i\over\tau}\right),\qquad i~=~1,2,\ldots,n-1.
 \end{array}
\eeq
Observe that $\delta_\ell(\xi)$ is supported in the dyadic cone
\bel{Cone}
\Lambda_\ell~=~ \left\{  (\tau,\lambda)\in\R\times\R^{n-1}~\colon~      2^{-\ell_i-1}~<~ {|\lambda_i|\over|\tau|}~<~2^{-\ell_i+1},~i=1,2,\ldots,n-1 \right\}.
\eeq
\begin{figure}[h]
\centering
\includegraphics[scale=0.12]{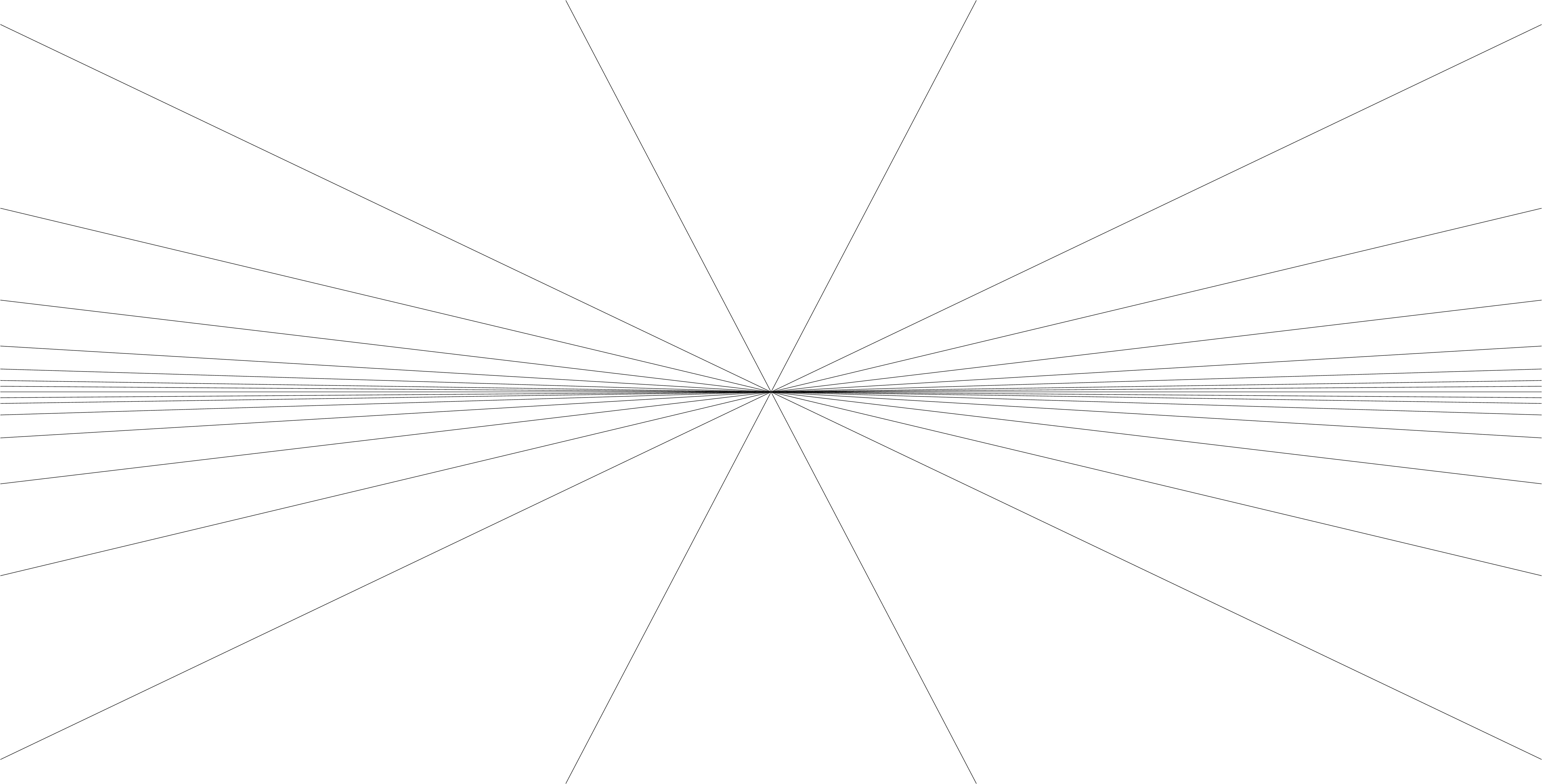}
\caption{$(\tau,\lambda)\in\R\times\R$ where $\tau$ is in the horizontal direction.}
\end{figure}

From (\ref{delta_t})-(\ref{Cone}) and direct computation, we find
\bel{delta Diff Ineq}
\left| \p_\tau^\alpha \p_\lambda^\beta  \delta(\tau,\lambda)\right|~\leq~\C_{\alpha~\beta}~ \left({1\over |\tau|}\right)^{\alpha}\prod_{i=1}^{n-1} \left({1\over |\lambda_i|}\right)^{\beta_i}
\eeq
for every multi-indices $\alpha, \beta$.

Define the partial operator
\bel{Partial}
\begin{array}{cc}\ds
\F_\ell f(x)~=~\int_{\R^n} f(y)\Omega_\ell(x,y)dy,
\\\\ \ds
\Omega_\ell(x,y)~=~\int_{\R^n}e^{2\pi\i \left(\Phi(x,\xi)-y\cdot\xi\right)}\sigma(x,\xi)\delta_\ell(\xi)d\xi.
\end{array}
\eeq

{\bf Lemma One}~~{\it 
 Suppose $\sigma\in\S^m$ for  $-n/2<m\leq0$. We have
\bel{F_t 2,p result}
\begin{array}{cc}\ds
 \left\| \F_\ell f\right\|_{\L^2(\R^n)}~\leq~\C_{p~\sigma~\Phi}~ \prod_{i=1}^{n-1}2^{\left({m\over n}\right) \ell_i}\left\| f\right\|_{\L^p(\R^n)}
\qquad
\hbox{for}\qquad{-m\over n}~=~{1\over p}-{1\over 2},
\end{array}
\eeq
\bel{F_t p',2 result}
\begin{array}{cc}\ds
 \left\| \F_\ell f\right\|_{\L^{p\over p-1}(\R^n)}~\leq~\C_{p~\sigma~\Phi}~\prod_{i=1}^{n-1}2^{\left({m\over n}\right) \ell_i}\left\| f\right\|_{\L^2(\R^n)}
\qquad
 \hbox{for}\qquad
{-m\over n}~=~{1\over 2}-{p-1\over p}.
\end{array}
\eeq}

In section 3, we prove {\bf Lemma One} together with  the $\L^2$-estimate:
\bel{L^2-result}
 \left\| \F f\right\|_{\L^2(\R^n)}~\leq~\C_{\sigma~\Phi}~\left\| f\right\|_{\L^2(\R^n)},\qquad \sigma\in\S^0. \eeq
Our main objective is to conclude 
\bel{H^1 L^1 F}
\begin{array}{lr}\ds
\left\|\F f\right\|_{\L^1(\R^n)}~\leq~\C_{\sigma~\Phi}~\left\|f\right\|_{\H^1(\R^n)},
\qquad
\left\|\F f\right\|_{\B\M\O(\R^n)}~\leq~\C_{\sigma~\Phi}~\left\|f\right\|_{\L^\infty(\R^n)}
\end{array}
\eeq
for $\sigma\in\S^{-{n-1\over 2}}$.

From (\ref{L^2-result}) and (\ref{H^1 L^1 F}), 
we can then  finish the proof of {\bf Theorem Two}  by carrying out an interpolation argument set out at {\bf 4.9}, chapter IX in the book of Stein \cite{Stein}.

Furthermore, by the duality between $\H^1$ and $\BMO$ spaces, as investigated by Fefferman \cite{Fefferman}, the second norm inequality in (\ref{H^1 L^1 F}) is equivalent to 
\bel{H^1 L^1 F*}
\left\|\F^* f\right\|_{\L^1(\R^n)}~\leq~\C_{\sigma~\Phi}~\left\|f\right\|_{\H^1(\R^n)},
\qquad \sigma\in\S^{-{n-1\over 2}}.
\eeq
Let $a$ to be an $\H^1$-{\it atom} associated to a ball $B_r(x_o)$ centered on some $x_o\in\R^n$ with  radius $r>0$. 
In order to obtain (\ref{H^1 L^1 F}), it is suffice to have
\bel{H^1 est F_ta}
\begin{array}{lr}\ds
\int_{\R^n} \left|\F a(x)\right|dx~\leq~\C_{\sigma~\Phi},\qquad \int_{\R^n} \left|\F^* a(x)\right|dx~\leq~\C_{\sigma~\Phi},\qquad\sigma\in\S^{-{n-1\over 2}}.
\end{array}
\eeq 
See 
 the characterization of $\H^1$-Hardy space established by Fefferman and Stein \cite{FC.S}.

Consider  a subset $\Q_r(x_o)\subset\R^n$, so-called the {\it region of influence},  satisfying
\bel{Q norm}
|\Q_r(x_o)|~\leq~\C_\sigma~r.
\eeq
Let $\F_\ell$ defined in (\ref{Partial}).
By using Schwartz inequality, we find
\bel{Local est}
\begin{array}{lr}\ds
\int_{\Q_r(x_o)} \left|\F_\ell a(x)\right|dx~\leq~|\Q_r(x_o)|^{1\over 2}\left\|\F_\ell a\right\|_{\L^2(\R^n)}
\\\\ \ds~~~~~~~~~~~~~~~~~~~~~~~~~~~~
~\leq~\C_\sigma~r^{1\over 2}\left\|\F_\ell a\right\|_{\L^2(\R^n)}.
\end{array}
\eeq
By applying {\bf Lemma One},  (\ref{F_t 2,p result}) implies
\bel{F_t a Result p 2}
\begin{array}{cc}\ds
\left\|\F_\ell a\right\|_{\L^2(\R^n)}~\leq~\C_{p~\sigma~\Phi}~\prod_{i=1}^{n-1} 2^{-\left({n-1\over 2n}\right)\ell_i}~ \|a\|_{\L^p(\R^n)},\qquad \sigma\in\S^{-{n-1\over 2}}
\\\\ \ds
\hbox{for}\qquad {1\over p}~=~{1\over 2}+{n-1\over 2n}.
\end{array}
\eeq
Note that 
$\|a\|_{\L^p(\R^n)} \leq |B_r(x_o)|^{-1+{1\over p}}$
because
$|a(x)|\leq |B_r(x_o)|^{-1}$ and $a$ is supported inside $B_r(x_o)$. Moreover,
$-1+{1\over p}=-{1\over 2}+{n-1\over 2n}=-{1\over 2n}$.

Together with (\ref{F_t a Result p 2}), we have
\bel{Local est 2}
\begin{array}{lr}\ds
\int_{\Q_r(x_o)} \left|\F_\ell a(x)\right|dx~\leq~\C_{p~\sigma~\Phi}~r^{1\over 2}\prod_{i=1}^{n-1}2^{-\left({n-1\over 2n}\right)\ell_i}~ \|a\|_{\L^p(\R^n)}
\\\\ \ds~~~~~~~~~~~~~~~~~~~~~~~~~~~
~\leq~\C_{p~\sigma~\Phi}~ r^{1\over 2}r^{n\left(-1+{1\over p}\right)}~\prod_{i=1}^{n-1} 2^{-\left({n-1\over 2n}\right)\ell_i}
~=~\C_{\sigma~\Phi}~\prod_{i=1}^{n-1}2^{-\left({n-1\over 2n}\right)\ell_i}.
\end{array}
\eeq
Clearly,  from (\ref{Local est 2}), we conclude 
\bel{Local est 3}
\begin{array}{lr}\ds
\int_{\Q_r(x_o)} \left|\F a(x)\right|dx
~\leq~\C_{\sigma~\Phi}~\sum_\ell~\prod_{i=1}^{n-1}2^{-\left({n-1\over 2n}\right)\ell_i}~\leq~\C_{\sigma~\Phi}.
\end{array}
\eeq
On the other hand,  by applying {\bf Lemma One}, (\ref{F_t p',2 result}) implies
\bel{F_t a Result 2 p' and F^*_t a Result p 2}
\begin{array}{cc}\ds
\left\|\F^*_\ell a\right\|_{\L^2(\R^n)}~\leq~\C_{p~\sigma~\Phi}~\prod_{i=1}^{n-1} 2^{-\left({n-1\over 2n}\right)\ell_i}~ \|a\|_{\L^p(\R^n)}
\\\\ \ds
\hbox{for}\qquad 
{1\over p}~=~{1\over 2}~+~{n-1\over 2n}.
\end{array}
\eeq
The {\it region of influence} associated to $\F^*$  is denoted by $\Q^*_r(x_o)$ satisfying
\bel{Q^* norm}
|\Q_r^*(x_o)|~\leq~\C_\sigma~r.
\eeq
By repeating the estimate in (\ref{Local est})-(\ref{Local est 3}) and using (\ref{F_t a Result 2 p' and F^*_t a Result p 2}) instead of (\ref{F_t a Result p 2}), we find
\bel{H^1 est *}
\int_{\Q^*_r(x_o)} \left|\F^* a(x)\right|dx~\leq~\C_{\sigma~\Phi}.
\eeq
Therefore, our task can be completed if we show
\bel{Comple est}
\begin{array}{lr}\ds
\int_{\R^n\setminus\Q_r(x_o)} \left|\F a(x)\right|dx~\leq~\C_{\sigma~\Phi},
\qquad
\int_{\R^n\setminus\Q^*_r(x_o)} \left|\F^* a(x)\right|dx~\leq~\C_{\sigma~\Phi}
\end{array}
\eeq
for $\sigma\in\S^{-{n-1\over 2}}$.

 In section 4, we give a heuristic estimate for (\ref{Comple est}) by assuming  that the kernel of the  partial operator  satisfies certain {\it majorization} properties. These properties are accumulated into {\bf Lemma Two}.

In section 5, we construct a second dyadic decomposition  in analogue to the framework of Seeger, Sogge and Stein \cite{S.S.S}.
The frequency space is asserted as an union of   geometric cones, denoted by $\Gamma_j^\nu$  whose central directions $\xi_j^\nu$ are almost uniformly distributed on  $\mathds{S}^{n-1}$ with a grid length approximately equal to $2^{-j/2}, j\ge0$.  In particular, we shall study the intersection $\Gamma_j^\nu\cap\Lambda_\ell\cap\{2^{j-1}\leq|\xi|\leq2^{j+1}\}$.

In section 6, we explicitly define $\Q_r(x_o)$  and $\Q^*_r(x_o)$. The size estimates in (\ref{Q norm}) and (\ref{Q^* norm}) hold respectively.

We prove {\bf Lemma Two} in the last section.

\section{Proof of Lemma One and the $\L^2$-boundedness of $\F$}
\setcounter{equation}{0}
First, we show that $\F$ of order zero is bounded on $\L^2(\R^n)$. 
By applying Plancherel theorem, our assertion  reduces to
\bel{Sf}
\mathcal{S} f(x)~=~\int_{\R^n} e^{2\pi\i\Phi(x,\xi)} \sigma(x,\xi) f(\xi) d\xi
\eeq 
whose adjoint operator is
\bel{S*f}
\mathcal{S}^* f(\xi)~=~\int_{\R^n} e^{-2\pi\i\Phi(x,\xi)} \bar{\sigma}(x,\xi) f(x) dx.
\eeq 
Let $\c$  be a small positive constant. We define an {\it narrow cone} as follows: suppose $\xi$ and $\eta$ belong to a same narrow cone and $|\eta|\leq|\xi|$. By writing $\eta=\rho\xi+\eta^\dagger$ for $0\leq\rho\leq1$ and $\eta^\dagger$ perpendicular to $\xi$, we require  $|\eta^\dagger|\leq\c\rho|\xi|$. The value of $\c$ depends on $\Phi$.

Clearly,  we can decompose the frequency space for which $\mathcal{S}$ or $\mathcal{S}^*$  can be written as a finite sum of  partial operators.  Each one of them has  a symbol  supported in such an narrow cone.  Note that our decomposition has no affection to the differentiation property $w.r.t~x$ of $\sigma(x,\xi)$.

Recall the estimate given at {\it 3.1.1}, chapter IX of Stein \cite{Stein}.
We have
\bel{Phi_x est}
\left|\nabla_x \Big(\Phi(x,\xi)-\Phi(x,\eta)\Big)\right|~\ge~\C_\Phi~|\xi-\eta|
\eeq
whenever $\xi$ and $\eta$ belong to a same narrow cone. 

Consider
\bel{S*Sf}
\mathcal{S}^*\mathcal{S} f(\xi)~=~\int_{\R^n} f(\eta)\mathfrak{S}^\sharp(\xi,\eta)d\eta
\eeq
where
\bel{S*S Kernel}
\mathfrak{S}^\sharp(\xi,\eta)~=~\int_{\R^n} e^{2\pi\i\left(\Phi(x,\eta)-\Phi(x,\xi)\right)} \sigma(x,\eta)\bar{\sigma}(x,\xi) dx.
\eeq
\begin{remark} In order to prove the $\L^2$-boundedness of $\mathcal{S}$ for $\sigma\in\S^0$, we assume that $\sigma(x,\xi)$ has a support inside an narrow cone defined as above. Moreover, it only satisfies the differential inequality $w.r.t~x$ inside (\ref{Class}). 
\end{remark}
Recall that $\sigma(x,\xi)$ has a compact support in $x$. Hence that $\mathfrak{S}^\sharp(\xi,\eta)$ is bounded in norm.

By using (\ref{Phi_x est}), 
an $N$-fold integration by parts  $w.r.t~x$ gives
\bel{Omega^sharp int by parts}
\begin{array}{lr}\ds
\left|\mathfrak{S}^\sharp(\xi,\eta)\right|~\leq~\C_{\Phi~N}~
\left|\xi-\eta\right|^{-N}
\left|\int_{\R^n} e^{2\pi\i\left(\Phi(x,\eta)-\Phi(x,\xi)\right)} \nabla_x^N\Big(\sigma(x,\eta)\bar{\sigma}(x,\xi)\Big)dx\right|
\end{array}
\eeq
for $\xi\neq\eta$. 

Together with (\ref{Omega^sharp int by parts}), we have
\bel{Omega^sharp est}
\left|\mathfrak{S}^\sharp(\xi,\eta)\right|~\leq~\C_{\sigma~\Phi~N}~\left({1\over 1+|\xi-\eta|}\right)^{N}
\eeq
for every $N\ge1$.

To conclude the $\L^2$-boundedness of $\mathcal{S}$, we write
\bel{L^2 est S*S} 
\begin{array}{lr}\ds
\left\|\mathcal{S}^*\mathcal{S} f\right\|_{\L^2(\R^n)}~=~\left\{\int_{\R^n}\left|\int_{\R^n} f(\eta)\mathfrak{S}^\sharp(\xi,\eta)d\eta\right|^2d\xi\right\}^{1\over 2}
\\\\ \ds~~~~~~~~~~~~~~~~~~~~~~
~=~\left\{\int_{\R^n}\left|\int_{\R^n} f(\xi-\zeta)\mathfrak{S}^\sharp(\xi,\xi-\zeta)d\zeta\right|^2d\xi\right\}^{1\over 2}\qquad (~ \zeta=\xi-\eta ~) 
\\\\ \ds~~~~~~~~~~~~~~~~~~~~~~
~\leq~\C~\int_{\R^n} \left\{\int_{\R^n} \left|f(\xi-\zeta)\right|^2 \Big|\mathfrak{S}^\sharp(\xi,\xi-\zeta)\Big|^2 d\xi \right\}^{1\over 2 }d\zeta\qquad \hbox{\small{by Minkowski integral inequality}}
\\\\ \ds~~~~~~~~~~~~~~~~~~~~~~
~\leq~\C_{\sigma~\Phi~N}~\int_{\R^n} \left\{\int_{\R^n} \left|f(\xi-\zeta)\right|^2 \left({1\over 1+|\zeta|}\right)^{2N} d\xi \right\}^{1\over 2 }d\zeta\qquad \hbox{\small{by (\ref{Omega^sharp est})}}
\\\\ \ds~~~~~~~~~~~~~~~~~~~~~~
~=~\C_{\sigma~\Phi~N}~\left\| f\right\|_{\L^2(\R^n)}\int_{\R^n}  \left({1\over 1+|\zeta|}\right)^{N}d\zeta
\\\\ \ds~~~~~~~~~~~~~~~~~~~~~~
~\leq~\C_{\sigma~\Phi}~\left\| f\right\|_{\L^2(\R^n)}\qquad \hbox{\small{for $N$  sufficiently large. }}
\end{array}
\eeq

Next, we begin to prove (\ref{F_t 2,p result}) in {\bf Lemma One}.

Suppose $\sigma\in\S^m$ for $-n/2<m<0$. Let $\delta_\ell(\xi)$ defined in (\ref{delta_t}).
We write
\bel{F_t decom}
\begin{array}{lr}\ds
\F_\ell f(x)~=~\int_{\R^n} e^{2\pi\i\Phi(x,\xi)}\sigma(x,\xi)\delta_\ell(\xi)\Hat{f}(\xi)d\xi
\\\\ \ds
~=~\int_{\R^n} e^{2\pi\i\Phi(x,\xi)}\sigma(x,\xi)\left(1+|\xi|^2\right)^{-{m\over2}} \left[\delta_\ell(\xi)\Hat{f}(\xi)\left(1+|\xi|^2\right)^{m\over2}\right]d\xi
\\\\ \ds
~=~ \prod_{i=1}^{n-1}2^{\left({m\over n}\right) \ell_i}\int_{\R^n} e^{2\pi\i\Phi(x,\xi)}\sigma(x,\xi)\left(1+|\xi|^2\right)^{-{m\over2}} \left\{\delta_\ell(\xi)\Hat{f}(\xi)\prod_{i=1}^{n-1}2^{-\left({m\over n}\right) \ell_i}\left(1+|\xi|^2\right)^{m\over2}\right\}d\xi
\\\\ \ds
~\doteq~\prod_{i=1}^{n-1}2^{\left({m\over n}\right) \ell_i}\int_{\R^n} e^{2\pi\i\Phi(x,\xi)}\sigma(x,\xi)\left(1+|\xi|^2\right)^{-{m\over2}} \Hat{T_\ell f}(\xi)d\xi.
\end{array}
\eeq
Observe that  $\sigma(x,\xi)(1+|\xi|^2)^{-{m\over 2}}\in\S^0$. By using the $\L^2$-boundedness of $\F$, it is suffice to prove (\ref{F_t 2,p result}) for
$T_\ell $ defined implicitly in (\ref{F_t decom}). 
 We have
\bel{convolution T_tf}
\begin{array}{cc}\ds
\Big(T_\ell f\Big)(x)~=~\int_{\R^n}f(y)\mathcal{K}_\ell(x-y)dy,
\\\\ \ds
\mathcal{K}_\ell(x)~=~\int_{\R^n} e^{2\pi\i x\cdot\xi} \delta_\ell(\xi) \prod_{i=1}^{n-1}2^{-\left({m\over n}\right)\ell_i}\left(1+|\xi|^2\right)^{m\over 2} d\xi.      
\end{array}     
\eeq
Let $\varphi$ be the smooth {\it bump}-function defined in (\ref{varphi}).  Consider
\bel{phi_j}
\phi_j(\xi)~=~\varphi\left(2^{-j}|\xi|\right)-\varphi\left(2^{-j+1}|\xi|\right),\qquad j\in\Z.
\eeq

Write $x=(z,w)\in\R\times\R^{n-1}$,  $y=(u,v)\in\R\times\R^{n-1}$ whose dual variable is $\xi=(\tau,\lambda)\in\R\times\R^{n-1}$. From (\ref{convolution T_tf})-(\ref{phi_j}), we have
\bel{K_t sum}
\begin{array}{lr}\ds
\mathcal{K}_\ell(x)
~=~ \prod_{i=1}^{n-1}2^{-\left({m\over n}\right) \ell_i}\int_{\R^n} e^{2\pi\i x\cdot\xi} \delta_\ell(\xi) \left(1+|\xi|^2\right)^{m\over 2} d\xi
\\\\ \ds~~~~~~~~~
~=~ \prod_{i=1}^{n-1}2^{-\left({m\over n}\right) \ell_i}\sum_{j\in\Z}\iint_{\R\times\R^{n-1}} e^{2\pi\i (z\tau+w\cdot\lambda)} \delta_\ell(\tau,\lambda) \phi_j(\tau,\lambda)\left(1+\tau^2+|\lambda|^2\right)^{m\over 2} d\tau d\lambda.
\end{array}
\eeq
Note that  $\phi_j(\xi)$ is supported in the dyadic annuli $2^{j-1}\leq|\xi|\leq2^{j+1}$. 
On the other hand, $\delta_\ell(\xi)$ defined in (\ref{delta_t}) is supported in the dyadic cone $\Lambda_\ell$  given in (\ref{Cone}).
We have $|\tau|\leq\C 2^j$ and $|\lambda_i|\leq\C 2^{j-\ell_i}, i=1,2,\ldots,n-1$ so that
\bel{support}
\left| \supp\delta_\ell(\tau,\lambda)\phi_j(\tau,\lambda)\right|~\leq~\C~ 2^j\prod_{i=1}^{n-1}2^{j-\ell_i}.
\eeq
Recall $\delta_\ell(\xi)$ satisfying  the differential inequality in (\ref{delta Diff Ineq}). We have
\bel{Computation Est1}
\begin{array}{lr}\ds
\left|\partial_\tau^M  \delta_\ell(\tau,\lambda) \phi_j(\tau,\lambda)\left(1+\tau^2+|\lambda|^2\right)^{m\over 2}\right|
~\leq~\C_M~\left(1+\tau^2+|\lambda|^2\right)^{m\over 2}\left({1\over |\tau|}\right)^M
\\\\ \ds~~~~~~~~~~~~~~~~~~~~~~~~~~~~~~~~~~~~~~~~~~~~~~~~~~~~~~~~~~~
~\leq~\C_M ~2^{jm}2^{-jM},
\\\\ \ds
\left|\partial_{\lambda_i}^{N_i}  \delta_\ell(\tau,\lambda) \phi_j(\tau,\lambda)\left(1+\tau^2+|\lambda|^2\right)^{m\over 2}\right|~\leq~\C_{N_i}~\left(1+\tau^2+|\lambda|^2\right)^{m\over 2}\left({1\over |\lambda_i|}\right)^{N_i}
\\\\ \ds~~~~~~~~~~~~~~~~~~~~~~~~~~~~~~~~~~~~~~~~~~~~~~~~~~~~~~~~~~~~
~\leq~\C_{N_i}~2^{jm}2^{-(j-\ell_i)N_i}
\end{array}
\eeq
for every $M\ge1$ and $N_i\ge1, i=1,2,\ldots,n-1$.

Let $N=N_1+N_2+\cdots+N_{n-1}$. An $M+N$-fold integration by parts $w.r.t~(\tau,\lambda)$ gives 
\bel{K_t by parts} 
\begin{array}{lr}\ds
\prod_{i=1}^{n-1}2^{-\left({m\over n}\right) \ell_i}\left|\iint_{\R\times\R^{n-1}} e^{2\pi\i (z\tau+w\cdot\lambda)} \delta_\ell(\tau,\lambda) \phi_j(\tau,\lambda)\left(1+\tau^2+|\lambda|^2\right)^{m\over 2} d\tau d\lambda\right|
\\\\ \ds
~\leq~\C_{M~N}~\prod_{i=1}^{n-1}2^{-\left({m\over n}\right) \ell_i} ~|z|^{-M}\prod_{i=1}^{n-1} |w_i|^{-N_i} 
\\\\ \ds~~~~~~~
\left| \iint_{\R\times\R^{n-1}} e^{2\pi\i (z\tau+w\cdot\lambda)} \partial_\tau^M\prod_{i=1}^{n-1}\p_{\lambda_i}^{N_i}\delta_\ell(\tau,\lambda) \phi_j(\tau,\lambda)\left(1+\tau^2+|\lambda|^2\right)^{m\over 2} d\tau d\lambda\right|
\\\\ \ds
~\leq~\C_{M~N}~ \prod_{i=1}^{n-1} 2^{-\left({m\over n}\right) \ell_i}\left\{2^{jm} 2^{j}\prod_{i=1}^{n-1}2^{j-\ell_i}\right\}\left(2^{j}|z|\right)^{-M}\prod_{i=1}^{n-1}\left(2^{j-\ell_i}|w_i|\right)^{-N_i} \qquad \hbox{\small{by (\ref{support})-(\ref{Computation Est1})}}
\\\\ \ds
~=~\C_{M~N}~2^{j\left({n+m\over n}\right)}  \left(2^{j}|z|\right)^{-M}\prod_{i=1}^{n-1} 2^{(j-\ell_i)\left({n+m\over n}\right)}\left(2^{j-\ell_i}|w_i|\right)^{-N_i}.
\end{array}
\eeq
We choose 
\bel{N,M tau}
\begin{array}{cc}\ds
M=0~~~\hbox{if}~~~|z|\leq2^{-j}\qquad \hbox{or}\qquad M=1 ~~~\hbox{if}~~~ |z|>2^{-j}; 
\\\\ \ds
N_i=0~~~ \hbox{if}~~~ |w_i|\leq2^{-j+\ell_i}\qquad \hbox{or}\qquad N_i=1~~~\hbox{if}~~~ |w_i|>2^{-j+\ell_i},\qquad i=1,2,\ldots,n-1.
\end{array}
\eeq
From (\ref{K_t sum}) and (\ref{K_t by parts}), we have
\bel{Sum K_t est} 
\begin{array}{lr}\ds
\left|\mathcal{K}_\ell(z,w)\right|~\leq~\C_{M~N} \sum_j 2^{j\left({n+m\over n}\right)}  \left(2^{j}|z|\right)^{-M}\prod_{i=1}^{n-1} 2^{(j-\ell_i)\left({n+m\over n}\right)}\left(2^{j-\ell_i}|w_i|\right)^{-N_i}
\\\\ \ds~~~~~~~~~~~~~~~~ 
~\leq~\C_{M~N}  \left\{\sum_j2^{j\left({n+m\over n}\right)}  \left(2^j|z|\right)^{-M}\right\}
\prod_{i=1}^{n-1} \left\{\sum_j 2^{(j-\ell_i)\left({n+m\over n}\right)}\left(2^{j-\ell_i}|w_i|\right)^{-N_i}\right\}
\\\\ \ds~~~~~~~~~~~~~~~~
~=~\C~  \left\{\sum_{|z|\leq2^{-j}} 2^{j\left({n+m\over n}\right)}~+~\sum_{|z|>2^{-j}}2^{j\left({n+m\over n}\right)} \left(2^{j}|z|\right)^{-1}\right\}
\\\\ \ds~~~~~~~~~~~~~~~~~~~~~
\prod_{i=1}^{n-1} \left\{\sum_{|w_i|\leq2^{-j+\ell_i}} 2^{(j-\ell_i)\left({n+m\over n}\right)}~+~\sum_{|w_i|>2^{-j+\ell_i}}2^{(j-\ell_i)\left({n+m\over n}\right)} \left(2^{j-\ell_i}|w_i|\right)^{-1}\right\}
\qquad\hbox{\small{by (\ref{N,M tau})}}
\\\\ \ds~~~~~~~~~~~~~~~~
~\leq~\C~ \left\{\left({1\over |z|}\right)^{n+m\over n}+\left({1\over |z|}\right) \sum_{|z|>2^{-j}} 2^{j\left({m\over n}\right)}\right\}\prod_{i=1}^{n-1}\left\{\left({1\over |w_i|}\right)^{\left({n+m\over n}\right)}+\left({1\over |w_i|}\right) \sum_{|w_i|>2^{-j+\ell_i}} 2^{(j-\ell_i)\left({m\over n}\right)}\right\}
\\ \ds~~~~~~~~~~~~~~~~~~~~~~~~~~~~~~~~~~~~~~~~~~~~~~~~~~~~~~~~~~~~~~~~~~~~~~~~~~~~~~~~~~~~~~~~~~~~~~~~~~~~~~~~~~~~~~~~~~~~~~~~~~~~~~~~~~~~
 \hbox{\small{( $m<0$ )}}
\\ \ds~~~~~~~~~~~~~~~~
~\leq~\C~ \left({1\over |z|}\right)^{n+m\over n}\prod_{i=1}^n\left({1\over|w_i|}\right)^{\left({n+m\over n}\right)}.
\end{array}
\eeq
Let $-{m\over n}={1\over p}-{1\over 2}$. By applying the Hardy-Littlewood-Sobolev inequality \cite{Hardy-Littlewood}-\cite{Sobolev} on every coordinate subspace and carrying out an iteration argument \footnote{See section 6 of \cite{Wang} for example.} by using Minkowski integral inequality,
we have
\bel{T_t regularity}
\begin{array}{lr}\ds
\left\|T_\ell f\right\|_{\L^2(\R^n)}~=~\left\{\iint_{\R\times\R^{n-1}} \Big(f\ast\mathcal{K}_\ell \Big)^2(z,w)dzdw\right\}^{1\over 2}
\\\\ \ds~~~~~~~~~~~~~~~~~~~
~\leq~\C~\left\{\iint_{\R\times\R^{n-1}} \left\{\iint_{\R\times\R^{n-1}} \left|f(u,v)\right|\left({1\over|z-u|}\right)^{n+m\over n}\prod_{i=1}^{n-1}\left({1\over|w_i-v_i|}\right)^{n+m\over n}dudv\right\}^2dzdw\right\}^{1\over 2}
\\ ~~~~~~~~~~~~~~~~~~~~~~~~~~~~~~~~~~~~~~~~~~~~~~~~~~~~~~~~~~~~~~~~~~~~~~~~~~~~~~~~~~~~~~~~~~~~~~~~~~~~~~~~~~~~~~~~~~~~~~~~~~~~~~~~~~~~~~~~~
 \hbox{\small{by (\ref{Sum K_t est})}}
\\ \ds~~~~~~~~~~~~~~~~~~
~\leq~\C_p~\left\| f\right\|_{\L^p(\R^n)}.
\end{array}
\eeq

Now, we turn to (\ref{F_t p',2 result}) in {\bf Lemma One}.

Consider
\bel{S_lf}
\mathcal{S}_\ell f(x)~=~\int_{\R^n} e^{2\pi\i\Phi(x,\xi)} \sigma(x,\xi)\delta_\ell(\xi) f(\xi) d\xi
\eeq 
where
\bel{S*_lf}
\mathcal{S}^* f(\xi)~=~\int_{\R^n} e^{-2\pi\i\Phi(x,\xi)} \bar{\sigma}(x,\xi)\bar{\delta}_\ell(\xi) f(x) dx.
\eeq 
Recall $\F_\ell$ defined in (\ref{Partial}). We have
\bel{F^*_t f}
\begin{array}{lr}\ds
\F^*_\ell f(x)~=~\int_{\R^n} f(y) \left\{\int_{\R^n} e^{2\pi\i \left(x\cdot\xi-\Phi(y,\xi)\right)}\bar{\sigma}(y,\xi)\bar{\delta}_\ell(\xi)d\xi\right\}dy
\\\\ \ds~~~~~~~~~~~
~=~\int_{\R^n} e^{2\pi\i x\cdot\xi} \mathcal{S}^*_\ell f(\xi)d\xi.
\end{array}
\eeq
We prove (\ref{F_t p',2 result})  by showing  
\bel{F_t p',2 result dual}
\begin{array}{cc}\ds
 \left\| \F^*_\ell f\right\|_{\L^2(\R^n)}~\leq~\C_{p~\sigma~\Phi}~\prod_{i=1}^{n-1}2^{\left({m\over n}\right) \ell_i}\left\| f\right\|_{\L^p(\R^n)}
\qquad
 \hbox{for}\qquad
{-m\over n}~=~{1\over p}-{1\over 2}.
\end{array}
\eeq
By using Plancherel theorem,  it is suffice to obtain (\ref{F_t p',2 result dual}) for $\mathcal{S}^*_\ell$.
Note that
 \bel{S^*_t L^2-norm}
\begin{array}{lr}\ds
\left\| \mathcal{S}^*_\ell f\right\|_{\L^2(\R^n)}^2~=~\int_{\R^n} \mathcal{S}_\ell\mathcal{S}^*_\ell f(x) f(x)dx
\\\\ \ds~~~~~~~~~~~~~~~~~~~~
~\leq~\left\|\mathcal{S}_\ell\mathcal{S}^*_\ell f\right\|_{\L^{p\over p-1}(\R^n)}\left\| f\right\|_{\L^p(\R^n)}\qquad\hbox{\small{by H\"{o}lder inequality.}}
\end{array}
\eeq
We aim to show
\bel{SS_t^* p',p result}
\begin{array}{cc}\ds
 \left\|\mathcal{S}_\ell \mathcal{S}^*_\ell f\right\|_{\L^{p\over p-1}(\R^n)}~\leq~\C_{p~\sigma~\Phi}~\prod_{i=1}^{n-1}2^{\left({2m\over n}\right) \ell_i}\left\| f\right\|_{\L^p(\R^n)}
 \\\\ \ds
\hbox{for}\qquad {-2m\over n}~=~{1\over p}-{p-1\over p}. 
\end{array}
\eeq
From (\ref{S_lf})-(\ref{S*_lf}), we have
\bel{S_tS^*_t}
\mathcal{S}_\ell\mathcal{S}^*_\ell f(x)~=~\int_{\R^n} f(y)\mathfrak{S}^\flat_\ell(x,y)dy
\eeq
where
\bel{SS* Kernel}
\mathfrak{S}^\flat_\ell(x,y)~=~\int_{\R^n} e^{2\pi\i\left(\Phi(x,\xi)-\Phi(y,\xi)\right)} \sigma(x,\xi)\delta_\ell(\xi)\bar{\sigma}(y,\xi)\bar{\delta}_\ell(\xi)d\xi.
\eeq
Write
\bel{nabla Phi}
\nabla_\xi \left(\Phi(x,\xi)-\Phi(y,\xi)\right)~=~\left[{\p^2\Phi\over \p x\p \xi}\right](x,\xi)(x-y)+\O\left(|x-y|^2\right).
\eeq
Note that $\Phi(x,\xi)$ satisfies the non-degeneracy condition  in (\ref{nondegeneracy}). From (\ref{nabla Phi}), we have
\bel{nabla Phi size}
\left|\nabla_\xi \left(\Phi(x,\xi)-\Phi(y,\xi)\right)\right|~\ge~\C_{\Phi}~|x-y|
\eeq
for $x,y$  sufficiently close.

\begin{remark} Recall that $\sigma(x,\xi)$ has a compact support in $x$.
By using a smooth partition of unity,  we can write it  as a finite sum of symbol functions. Each one of them has a sufficiently small $x$-support  depending on $\Phi$. 
 We make this assumption for the remaining section. Note that it does not affect the differentiation $w.r.t~\xi$ for $\sigma(x,\xi)$.
\end{remark}
Consider 
\bel{t dila}
\begin{array}{cc}\ds
x~=~\mathfrak{L}^{-1}x'~=~\left(z, 2^{\ell_1}w'_1,\ldots,2^{\ell_{n-1}}w'_{n-1}\right),\qquad y~=~\mathfrak{L}^{-1}y'~=~\left(u,2^{\ell_1}v'_1,\ldots,2^{\ell_{n-1}}v'_{n-1}\right)
\\\\ \ds
\hbox{and}\qquad \xi~=~\mathfrak{L}\xi'~=~\left(\tau, 2^{-\ell_1}\lambda'_1,\ldots,2^{-\ell_{n-1}}\lambda'_{n-1}\right).
\end{array}
\eeq
By definition of $\delta_\ell(\xi)$ in (\ref{delta_t}), we have
$\delta_\ell(\mathfrak{L}\xi')=\delta_o(\xi')$.
Let $\phi_j(\xi)$ defined in (\ref{phi_j}). We write
\bel{Sum S^flat dila}
\begin{array}{lr}\ds
\mathfrak{S}^\flat_\ell\left(x,y\right)~=~\mathfrak{S}^\flat_\ell\left(\mathfrak{L}^{-1}x',\mathfrak{L}^{-1}y'\right)
\\\\ \ds
~=~\prod_{i=1}^{n-1}2^{-\ell_i}\sum_{j\in\Z} \int_{\R^n} e^{2\pi\i\left(\Phi(\mathfrak{L}^{-1}x',\mathfrak{L}\xi')-\Phi(\mathfrak{L}^{-1}y',\mathfrak{L}\xi')\right)} \sigma(\mathfrak{L}^{-1}x',\mathfrak{L}\xi')\delta_o(\xi')\bar{\sigma}(\mathfrak{L}^{-1}y',\mathfrak{L}\xi')\bar{\delta}_o(\xi')\phi_j(\xi')d\xi'.
\end{array}
\eeq
Note that  $2^{j-1}\leq|\xi'|\leq2^{j+1}$ for $\xi'$  in the support of $\phi_j(\xi')$. Recall $\sigma\in\S^m$ satisfying the differential inequality in (\ref{Class}). We have
\bel{Computation Est2}
\begin{array}{rl}
\left|\p^\alpha_{\xi'} \sigma\left(\mathfrak{L}^{-1}x',\mathfrak{L}\xi'\right)\delta_o(\xi')\bar{\sigma}\left(\mathfrak{L}^{-1}y',\mathfrak{L}\xi'\right)\bar{\delta}_o(\xi')\phi_j(\xi')\right|
~\leq~\C_\alpha~2^{-j|\alpha|} 2^{j m}
\end{array}
\eeq
for every multi-index $\alpha$.

Moreover, 
\bel{det Phi}
\det\left[{\p^2\Phi\over \p x\p \xi}\right](x,\xi)~=~\det\left[{\p^2\Phi\over \p x'\p \xi'}\right](\mathfrak{L}^{-1}x',\mathfrak{L}\xi').
\eeq
Indeed,    $2^{\ell_i}$ appears at  the $i$-th column of $\left[{\p^2\Phi\over \p x'\p \xi'}\right](\mathfrak{L}^{-1}x',\mathfrak{L}\xi')$. On the other hand, $2^{-\ell_i}$ appears at  every $i$-th row of $\left[{\p^2\Phi\over \p x'\p \xi'}\right](\mathfrak{L}^{-1}x',\mathfrak{L}\xi')$ respectively. 

From (\ref{nabla Phi})-(\ref{nabla Phi size}) and (\ref{det Phi}), we find
\bel{nabla Phi size new}
\left|\nabla_{\xi'} \left(\Phi(\mathfrak{L}^{-1}x',\mathfrak{L}\xi')-\Phi(\mathfrak{L}^{-1}y',\mathfrak{L}\xi')\right)\right|~\ge~\C_{\Phi}~|x'-y'|.
\eeq
By using (\ref{Computation Est2}) and (\ref{nabla Phi size new}), an $M+N$-fold 
integration by parts $w.r.t~ \xi'$ shows
\bel{S^flat by parts}
\begin{array}{lr}\ds
 \prod_{i=1}^{n-1}2^{-\ell_i}\left|\int_{\R^n} e^{2\pi\i\left(\Phi(\mathfrak{L}^{-1}x',\mathfrak{L}\xi')-\Phi(\mathfrak{L}^{-1}y',\mathfrak{L}\xi')\right)} \sigma(\mathfrak{L}^{-1}x',\mathfrak{L}\xi')\delta_o(\xi')\bar{\sigma}(\mathfrak{L}^{-1}y',\mathfrak{L}\xi')\bar{\delta}_o(\xi')\phi_j(\xi')d\xi'\right|
\\\\ \ds
~\leq~\C_{\Phi~M~N}~ \prod_{i=1}^{n-1}2^{-\ell_i}\left(2^{2mj} 2^{jn}\right)\left(2^j |x'-y'|\right)^{-M-N}
\\\\ \ds
~\leq~\C_{\Phi~M~N}~\left( 2^{2mj}\right)\left\{ 2^{j}\prod_{i=1}^{n-1}2^{j-\ell_i}\right\} \left( 2^j |z-u|\right)^{-M}\prod_{i=1}^{n-1}\left( 2^j |w'_i-v'_i|\right)^{-N_i}
\\\\ \ds
~=~\C_{\Phi~M~N}~ \prod_{i=1}^{n-1}2^{\left({2m\over n}\right)\ell_i}\left[2^{j\left({n+2m\over n}\right)}  \left( 2^j |z-u|\right)^{-M}\right]\prod_{i=1}^{n-1}\left[2^{(j-\ell_i)\left({n+2m\over n}\right)}  \left( 2^{j-\ell_i} |w_i-v_i|\right)^{-N_i}\right].
\end{array}
\eeq
We choose 
\bel{N,M lambda}
\begin{array}{cc}\ds
M=0~~~\hbox{if}~~~ |z-u|\leq2^{-j}\qquad \hbox{or}\qquad N=1~~~\hbox{if}~~~ |z-u|>2^{-j},
\\\\ \ds
N_i=0~~~\hbox{if}~~~ |w_i-v_i|\leq2^{-j+\ell}\qquad \hbox{or}\qquad N_i=1~~~\hbox{if}~~~ |w_i-v_i|>2^{-j+\ell_i},\qquad i=1,2,\ldots,n-1.
\end{array}
\eeq
From (\ref{Sum S^flat dila}) and (\ref{S^flat by parts}), we have
\bel{Sum S^flat est} 
\begin{array}{lr}\ds
\left|\mathfrak{S}^\flat_\ell\left(x,y\right)\right|~\leq~\C_{\Phi~M~N}~ \sum_{j\in\Z} ~\prod_{i=1}^{n-1}2^{\left({2m\over n}\right)\ell_i}\left\{2^{j\left({n+2m\over n}\right)}  \left( 2^j |z-u|\right)^{-M}\right\}\prod_{i=1}^{n-1}\left\{2^{(j-\ell_i)\left({n+2m\over n}\right)}  \left( 2^{j-\ell_i} |w_i-v_i|\right)^{-N_i}\right\}
\\\\ \ds~~~~~~~~~~~~~~~ 
~\leq~\C_{\Phi~M~N}~  \prod_{i=1}^{n-1}2^{\left({2m\over n}\right)\ell_i} \left\{ \sum_{j\in\Z} 2^{j\left({n+2m\over n}\right)}  \left( 2^{j} |z-u|\right)^{-M} \right\}\prod_{i=1}^{n-1}\left\{\sum_{j\in\Z} 2^{(j-\ell_i)\left({n+2m\over n}\right)}  \left( 2^{j-\ell_i} |w_i-v_i|\right)^{-N_i}\right\}
\\\\ \ds~~~~~~~~~~~~~~~
~=~\C_{\Phi}~ \prod_{i=1}^{n-1}2^{\left({2m\over n}\right)\ell_i} \left\{\sum_{|z-u|\leq2^{-j}} 2^{j\left({n+2m\over n}\right)}~+~\sum_{|z-u|>2^{-j}}2^{j\left({n+2m\over n}\right)} \left(2^{j}|z-u|\right)^{-1}\right\}
\\\\ \ds~~~~~~~~~~~~~~~~~~~~~~~~~~~~
\prod_{i=1}^{n-1}\left\{\sum_{|w_i-v_i|\leq2^{-j+\ell_i}} 2^{(j-\ell_i)\left({n+2m\over n}\right)}~+~\sum_{|w_i-v_i|>2^{-j+\ell_i}} 2^{(j-\ell_i)\left({n+2m\over n}\right)} \left(2^{j-\ell_i}|w_i-v_i|\right)^{-1}\right\}~~\hbox{\small{by (\ref{N,M lambda})}}
\\\\ \ds~~~~~~~~~~~~~~~
~\leq~\C_\Phi~ \prod_{i=1}^{n-1}2^{\left({2m\over n}\right)\ell_i} \left\{ \left({1\over |z-u|}\right)^{n+2m\over n}~+~\left({1\over |z-u|}\right) \sum_{|z-u|>2^{-j}} 2^{j\left({2m\over n}\right)}\right\}
\\\\ \ds~~~~~~~~~~~~~~~~~~~~~~~~~~~~
 \prod_{i=1}^{n-1}\left\{ \left({1\over |w_i-v_i|}\right)^{{n+2m\over n}}~+~\left({1\over |w_i-v_i|}\right) \sum_{|w_i-v_i|>2^{-j+\ell_i}} 2^{(j-\ell_i)\left({2m\over n}\right)}\right\}
\qquad \hbox{\small{($m<0$)}}
\\\\ \ds~~~~~~~~~~~~~~~
~\leq~\C_\Phi~ \prod_{i=1}^{n-1}2^{\left({2m\over n}\right)\ell_i}~\left({1\over |z-u|}\right)^{n+2m\over n}\prod_{i=1}^{n-1}\left({1\over |w_i-v_i|}\right)^{{n+2m\over n}}.
\end{array}
\eeq
Let $ -{2m\over n}={1\over p}-{p-1\over p}$. By applying Hardy-Littlewood-Sobolev inequality \cite{Hardy-Littlewood}-\cite{Sobolev} on every coordinate subspace and using Minkowski integral inequality, we have 
\bel{S_tS^*_t regularity}
\begin{array}{lr}\ds
\left\|\mathcal{S}_\ell\mathcal{S}^*_\ell f\right\|_{\L^{p\over p-1}(\R^n)}~=~\left\{\int_{\R^n}\left|\int_{\R^n} f(y)\mathfrak{S}^\flat_\ell(x,y)dy\right|^{p\over p-1}dx\right\}^{p-1\over p}\qquad\hbox{\small{by (\ref{S_tS^*_t})}}
\\\\ \ds
~\leq~\C_\Phi \prod_{i=1}^{n-1}2^{\left({2m\over n}\right)\ell_i}~\left\{\iint_{\R\times\R^{n-1}}\left|\iint_{\R\times\R^{n-1}} |f(u,v)|\left({1\over |z-u|}\right)^{n+2m\over n}\prod_{i=1}^{n-1}\left({1\over |w_i-v_i|}\right)^{{n+2m\over n}}dudv\right|^{p\over p-1}dzdw\right\}^{p-1\over p}
\\ \ds~~~~~~~~~~~~~~~~~~~~~~~~~~~~~~~~~~~~~~~~~~~~~~~~~~~~~~~~~~~~~~~~~~~~~~~~~~~~~~~~~~~~~~~~~~~~~~~~~~~~~~~~~~~~~~~~~~~~~~~~~~~~~~~\hbox{\small{by (\ref{Sum S^flat est})}}
\\ \ds
~\leq~\C_{p~\Phi}~ \prod_{i=1}^{n-1} 2^{\left({2m\over n}\right)\ell_i}~\left\| f\right\|_{\L^p(\R^n)}.
\end{array}
\eeq
Recall {\bf Remark 3.2}.  By using Minkowski inequality and (\ref{S_tS^*_t regularity}), we obtain (\ref{SS_t^* p',p result}) as desired.

\section{A heuristic estimate}
\setcounter{equation}{0}
Let $\I\cup\J=\{1,2,\ldots,n-1\}$ such that
\bel{IJ}
\begin{array}{cc}
  0~\leq~\ell_i~\leq~ j/2+3,\qquad i\in\I,
\qquad
  \ell_i~>~j/2+3,\qquad i\in\J,
\\\\ \ds
\J^\sharp~=~\{ i\in\J~\colon~\ell_i>j+3\},\qquad \J^\flat~=~\left\{ i\in\J~\colon~j/2+3<\ell_i\leq j+3\right\}
\end{array}
\eeq
for every $j>0$ and $\ell_i\ge0, i=1,2,\ldots,n-1$.

Their cardinalities are denoted by $|\I|$, $|\J|$, $|\J^\flat|$ and $|\J^\sharp|$ respectively.

Let $\varphi$ be the smooth {\it bump}-function  given in (\ref{varphi}). Recall
$\delta_\ell(\xi)$ defined in (\ref{delta_t})-(\ref{Cone}) and $\phi_j(\xi)$ defined in (\ref{phi_j}). 
Note that $\sum_\ell \delta_\ell(\xi)=\sum_{j}\phi_j(\xi)\equiv1$.

Define
\bel{delta_lj}
\begin{array}{lr}\ds
\delta_{\ell j}(\xi)~=~\prod_{i\in\I\cup\J^\flat}\delta_{\ell_i}(\xi)  \prod_{i\in\J^\sharp}~ \sum_{\ell_i} \delta_{\ell_i}(\xi)
\\\\ \ds~~~~~~~~~
~=~\prod_{i\in\I\cup\J^\flat}\delta_{\ell_i}(\xi)\prod_{i\in\J^\sharp}~ \sum_{\ell_i} \varphi\left(2^{\ell_i} {\lambda_i\over\tau}\right)-\varphi\left(2^{\ell_i+1}{\lambda_i\over \tau}\right)
\\\\ \ds~~~~~~~~~
~=\left.\begin{array}{lr}\ds  \prod_{i\in\I\cup\J^\flat}\delta_{\ell_i}(\xi)\prod_{i\in\J^\sharp}\varphi\left(2^{j+4} {\lambda_i\over\tau}\right) \qquad j+4>0,
\\\\ \ds
 \prod_{i\in\I\cup\J^\flat}\delta_{\ell_i}(\xi) \prod_{i\in\J^\sharp}\varphi\left( {\lambda_i\over\tau}\right)\qquad~~~~~~ j+4\leq0.
\end{array}\right.
\end{array}
\eeq
Let $j>0$. Observe that $\delta_{\ell j}(\xi)$ is supported in
\bel{Lambda_lj}
\begin{array}{lr}\ds
\Lambda_{\ell j}~=~\Bigg\{  (\tau,\lambda)\in\R\times\R^{n-1}~\colon~      2^{-\ell_i-1}~<~ {|\lambda_i|\over|\tau|}~<~2^{-\ell_i+1},~i\in\I\cup\J^\flat 
\\\\ \ds~~~~~~~~~~~~~~~~~~~~~~~~~~~~~~~~
\hbox{and}\qquad 0~<~ {|\lambda_i|\over|\tau|}~<~2^{-j-3},~i\in\J^\sharp\Bigg\}.
\end{array}
\eeq
Moreover, by definition of $\delta_\ell(\xi)$ in (\ref{delta_t}) and $\delta_{\ell j}(\xi)$ in (\ref{delta_lj}), we have
\bel{delta Sum}
\sum_{\ell_i\leq j+3, ~i=1,2,\ldots,n-1} \delta_{\ell j}(\xi)~\equiv~1.
\eeq
Consider
\bel{Omega_t,j}
\Omega_{\ell j}(x,y)~=~\int_{\R^n} e^{2\pi\i\left(\Phi(x,\xi)-y\cdot\xi\right)}\delta_{\ell j}(\xi)\phi_j(\xi)\sigma(x,\xi)d\xi
\eeq
where $\phi_j(\xi), j\in\Z$ is defined in (\ref{phi_j}).

Let $a$ be an $\H^1$-{\it atom} associated to the ball $B_r(x_o)$.
Recall $\F$ defined in (\ref{Ff}). We have
\bel{F_t f rewrite j}
\begin{array}{lr}\ds
\int_{\R^n\setminus\Q_r(x_o)}\left|\F a(x)\right| dx
~=~\int_{\R^n\setminus\Q_r(x_o)}\left|\int_{\R^n}a(y)\sum_{j}\sum_{\ell_i\leq j+3, ~i=1,2,\ldots,n-1}\Omega_{\ell j}(x,y)dy\right| dx
\\ \ds~~~~~~~~~~~~~~~~~~~~~~~~~~~~~~~~~~~~~~~~~~~~~~~~~~~~~~~~~~~~~~~~~~~~~~~~~~
~~~~~~~~~~~~~~~~~~~~
 \hbox{\small{by (\ref{delta_lj}) and (\ref{delta Sum})}}
\\\\ \ds~~~~~~~~~~~~~~~~~~~~~~~~~~~~~~~
~\leq~\int_{\supp\sigma}\left\{\int_{\R^n} |a(y)|\Bigg|\sum_{j\leq0}\sum_{\ell_i\leq j+3, ~i=1,2,\ldots,n-1}\Omega_{\ell j}(x,y)\Bigg|dy\right\}dx
\\\\ \ds~~~~~~~~~~~~~~~~~~~~~~~~~~~~~~~
~+~\sum_{j>0}\sum_{\ell_i\leq j+3, ~i=1,2,\ldots,n-1}~\int_{\R^n\setminus\Q_r(x_o)}\left|\int_{\R^n} a(y)\Omega_{\ell j}(x,y)dy\right|dx.
\end{array}
\eeq
Note that $\sum_{j\leq0}\phi_j(\xi)$ is supported inside the ball $|\xi|\leq2$.  
From (\ref{delta Sum})-(\ref{Omega_t,j}), we find 
the first term on the R.H.S of (\ref{F_t f rewrite j})  bounded by $\C_\sigma$.
\v

{\bf Lemma Two}~~{\it  
Suppose $\sigma\in\S^{-{n-1\over 2}}$. For every $ j>0$, we have
\bel{Est1}
\int_{\R^n} \left|\Omega_{\ell j}(x,y)\right| dx~\leq~\C_{\sigma~\Phi}~
\prod_{i\in\I}2^{-\ell_i}\prod_{j\in\J^\flat} 2^{-\left({1\over 2}\right)\ell_i}\prod_{j\in\J^\sharp} 2^{-\left({1\over 2}\right)j},
\eeq
\bel{Est2}
\int_{\R^n}\left|\Omega_{\ell j}(x,y)-\Omega_{\ell j}(x,x_o)\right|dx~\leq~\C_{\sigma~\Phi}~2^j|y-x_o|~
\prod_{i\in\I}2^{-\ell_i}\prod_{j\in\J^\flat} 2^{-\left({1\over 2}\right)\ell_i}\prod_{j\in\J^\sharp} 2^{-\left({1\over 2}\right)j}
\eeq
and
\bel{Est3}
\int_{\R^n\setminus\Q_r(x_o)}\left|\Omega_{\ell j}(x,y)\right|dx~\leq~\C_{\sigma~\Phi}~{2^{-j}\over r}~\prod_{i\in\I}2^{-\ell_i}\prod_{j\in\J^\flat} 2^{-\left({1\over 2}\right)\ell_i}\prod_{j\in\J^\sharp} 2^{-\left({1\over 2}\right)j},\qquad y\in B_r(x_o)
\eeq
whenever $2^{j}>r^{-1}$.}
\v
Consider $2^j\leq r^{-1}$.  We  write
\bel{F_t cancella}
\int_{\R^n} a(y)\Omega_{\ell j}(x,y)dy~=~\int_{B_r(x_o)} a(y) \left(\Omega_{\ell j}(x,y)-\Omega_{\ell j}(x,x_o)\right)dy
\eeq 
because 
$\int_{B_r(x_o)} a(y)dy=0$ and $a$ is supported in $B_r(x_o)$. 

By using (\ref{Est2}) and (\ref{F_t cancella}), we find
\bel{Norm Est1}
\begin{array}{lr}\ds
\int_{\R^n}\left|\int_{\R^n} a(y)\Omega_{\ell j}(x,y)dy\right|dx
~\leq~\int_{B_r(x_o)} |a(y)|\left\{\int_{\R^n} \left|\Omega_{\ell j}(x,y)-\Omega_{\ell j}(x,x_o)\right| dx\right\} dy
\\\\ \ds~~~~~~~~~~~~~~~~~~~~~~~~~~~~~~~~~~~~~~~~~~~~~~
~\leq~ \C_{\sigma~\Phi}~2^j|y-x_o|~\prod_{i\in\I}2^{-\ell_i}\prod_{j\in\J^\flat} 2^{-\left({1\over 2}\right)\ell_i}\prod_{j\in\J^\sharp} 2^{-\left({1\over 2}\right)j}
\\\\ \ds~~~~~~~~~~~~~~~~~~~~~~~~~~~~~~~~~~~~~~~~~~~~~~
~\leq~\C_{\sigma~\Phi}~2^j r~\prod_{i\in\I}2^{-\ell_i}\prod_{j\in\J^\flat} 2^{-\left({1\over 2}\right)\ell_i}\prod_{j\in\J^\sharp} 2^{-\left({1\over 2}\right)j}\qquad y\in B_r(x_o).
\end{array}
\eeq

By summing over all regarding $\ell$ and $j$ s, we have
\bel{Sum1}
\begin{array}{lr}\ds
 \sum_{2^j\leq r^{-1}} \sum_{\ell_i\leq j+3, ~i=1,2,\ldots,n-1}\int_{\R^n}\left|\int_{\R^n} a(y)\Omega_{\ell j}(x,y)dy\right|dx
\\\\ \ds
~\leq~\C_{\sigma~\Phi}~r\sum_{2^j\leq r^{-1}} 2^j 

\sum_{\I\cup\J^\flat\cup\J^\sharp=\{1,2,\ldots,n-1\}}
\left\{\sum_{\ell_i\leq j+3, ~i=1,2,\ldots,n-1} \prod_{i\in\I}2^{-\ell_i}\prod_{j\in\J^\flat} 2^{-\left({1\over 2}\right)\ell_i}\prod_{j\in\J^\sharp} 2^{-\left({1\over 2}\right)j}\right\}
\\ \ds~~~~~~~~~~~~~~~~~~~~~~~~~~~~~~~~~~~~~~~~~~~~~~~~~~~~~~~~~~~~~~~~~
~~~~~~~~~~~~~~~~~~~~~~~~~~~~~~~~~~~~~~~~~~
\hbox{\small{by (\ref{F_t cancella})-(\ref{Norm Est1})}}
\\ \ds
~\leq~\C_{\sigma~\Phi} ~r\sum_{2^j\leq r^{-1}} 2^j~\sum_{\I\cup\J^\flat\cup\J^\sharp=\{1,2,\ldots,n-1\}} \prod_{i\in\J^\sharp}j 2^{-\left({1\over 2}\right)j}
\\\\ \ds
~\leq~\C_{\sigma~\Phi} ~r\sum_{2^j\leq r^{-1}} 2^j~\leq~\C_{\sigma~\Phi}.
\end{array}
\eeq
For $2^j>r^{-1}$,  (\ref{Est3}) implies
\bel{Norm Est2}
\begin{array}{lr}\ds
\int_{\R^n\setminus\Q_r(x_o)}\left|\int_{\R^n} a(y)\Omega_{\ell j}(x,y)dy\right|dx
~\leq~\int_{B_r(x_o)} |a(y)|\left\{\int_{\R^n\setminus\Q_r(x_o)} \left|\Omega_{\ell j}(x,y)\right| dx\right\} dy
\\\\ \ds~~~~~~~~~~~~~~~~~~~~~~~~~~~~~~~~~~~~~~~~~~~~~~~~~~~~~~~
~\leq~\C_{\sigma~\Phi}~{2^{-j}\over r}~\prod_{i\in\I}2^{-\ell_i}\prod_{j\in\J^\flat} 2^{-\left({1\over 2}\right)\ell_i}\prod_{i\in\J^\sharp}2^{-\left({1\over 2}\right)j}.
\end{array}
\eeq
By summing over all regarding $\ell$ and $j$ s, we have
\bel{Sum2}
\begin{array}{lr}\ds
\sum_{2^j>r^{-1}} \sum_{\ell_i\leq j+3, ~i=1,2,\ldots,n-1}\int_{\R^n\setminus\Q_r(x_o)}\left|\int_{\R^n} a(y)\Omega_{\ell j}(x,y)dy\right|dx
\\\\ \ds
~\leq~\C_{\sigma~\Phi}~r^{-1}\sum_{2^j>r^{-1}} 2^{-j}~\sum_{\I\cup\J^\flat\cup\J^\sharp=\{1,2,\ldots,n-1\}}
\left\{\sum_{\ell_i\leq j+3, ~i=1,2,\ldots,n-1} \prod_{i\in\I}2^{-\ell_i}\prod_{j\in\J^\flat} 2^{-\left({1\over 2}\right)\ell_i}\prod_{j\in\J^\sharp} 2^{-\left({1\over 2}\right)j}\right\}
\\ \ds~~~~~~~~~~~~~~~~~~~~~~~~~~~~~~~~~~~~~~~~~~~~~~~~~~~~~~~~~~~~~~~~~
~~~~~~~~~~~~~~~~~~~~~~~~~~~~~~~~~~~~~~~~~~~~~~~~~~~
\hbox{\small{by (\ref{Norm Est2})}}
\\\\ \ds
~\leq~\C_{\sigma~\Phi}~r^{-1}\sum_{2^j>r^{-1}} 2^{-j}~\sum_{\I\cup\J^\flat\cup\J^\sharp=\{1,2,\ldots,n-1\}} \prod_{i\in\J^\sharp}j 2^{-\left({1\over 2}\right)j}
\\\\ \ds
~\leq~\C_{\sigma~\Phi}~r^{-1}\sum_{2^j>r^{-1}} 2^{-j}~\leq~\C_{\sigma~\Phi}.
\end{array}
\eeq
From (\ref{F_t f rewrite j}), 
(\ref{Sum1}) and (\ref{Sum2}), we obtain the first inequality in (\ref{Comple est}). 

On the other hand, define
\bel{Omega^*_t,j}
\Omega^*_{\ell j}(x,y)~=~\int_{\R^n} e^{2\pi\i\left(x\cdot\xi-\Phi(y,\xi)\right)}\bar{\delta}_{\ell j}(\xi)\bar{\phi}_j(\xi)\bar{\sigma}(y,\xi)d\xi
\eeq
for the associated adjoint operator. 
\begin{remark} $\Omega^*_{\ell j}(x,y)$  satisfies (\ref{Est1})-(\ref{Est3}) with $\Q_r(x_o)$ replaced by $\Q^*_r(x_o)$.
\end{remark}
We prove the second inequality in (\ref{Comple est})
 by  repeating the estimate in (\ref{F_t f rewrite j})-(\ref{Sum2}) with $\Omega_{\ell j}(x,y)$ and $\Q_r(x_o)$ replaced by $\Omega^*_{\ell j}(x,y)$ and $\Q^*_r(x_o)$.

\section{A second dyadic decomposition}
\setcounter{equation}{0}
For $\xi=(\tau,\lambda)\in\R\times\R^{n-1}$,  we denote 
\bel{Intersection k}
\begin{array}{cc}\ds
\mathbb{S}^{n-2}_n~\doteq~\mathbb{S}^{n-1}\cap\left\{(\tau,\lambda)\in\R\times\R^{n-1}~\colon~\tau=0\right\},
\\\\ \ds
\mathbb{S}^{n-2}_i~\doteq~\mathbb{S}^{n-1}\cap\left\{(\tau,\lambda)\in\R\times\R^{n-1}~\colon~\lambda_i=0\right\},\qquad i=1,2,\ldots,n-1.
 \end{array}
\eeq
Let $j>0$  fixed.   
We construct a collection of points $\{\xi^\nu_j\}_\nu\subset\mathds{S}^{n-1}$ as follows.

{\bf ( 1 )} Every unit vector $(\xi_i, \xi_i^\dagger)=(\pm1,0)\in\R\times\R^{n-1}, i=1,2,\ldots,n$ belongs to  $\{\xi^\nu_j\}_\nu$.

{\bf ( 2 )} For each $i=1,2,\ldots,n$,  a subset of  $\{\xi^\nu_j\}_\nu$ are equally distributed on $\mathbb{S}^{n-2}_i\subset\R^{n-1}$ with a grid length equal to $\c 2^{-j/2}$ for  $1/\sqrt{2}\leq\c\leq\sqrt{2}$. 

{\bf ( 3 )} The remaining of $\{\xi^\nu_j\}_\nu$ are equally distributed on $\mathds{S}^{n-1}\setminus\Cup_{i=1}^n \mathds{S}^{n-2}_i$ with the same grid length. 
\begin{remark}
From {\bf ( 1 )-( 3 )}, there are at most a constant multiple of $2^{j\left({n-1\over2}\right)}$ elements in $\{\xi^\nu_j\}_\nu$.
\end{remark}
 \begin{remark}
For every given $\xi\in\R^n$, there exists a $\xi^\nu_j$  such that $\left|{\xi\over|\xi|}-\xi^\nu_j\right|\leq2^{-j/2}$.
\end{remark}
Define 
\bel{Gamma}
\Gamma_{j}^{\nu}~=~\Bigg\{\xi\in\R^n~\colon~\left| {\xi\over |\xi|}-\xi^{\nu}_{j}\right|~\leq~3\cdot2^{-j/2}\Bigg\}
\eeq
whose central direction is $\xi^{\nu}_{j}$.
We have
\bel{intersection of cone norm}
\begin{array}{lr}\ds
\left|~\Gamma^\nu_j\cap\left\{2^{j-1}\leq|\xi|<2^{j+1}\right\}~\right|
~\leq~\C~2^j2^{j\left({n-1\over 2}\right)}.
\end{array}
\eeq
Recall $\I\cup\J^\flat\cup\J^\sharp=\{1,2,\ldots,n-1\}$ from (\ref{IJ}). Let $\Lambda_{\ell j}$ defined in (\ref{Lambda_lj}). We have
\bel{Lambda norm><}
\left|~\Lambda_{\ell j}\cap\left\{2^{j-1}\leq|\xi|<2^{j+1}\right\}~\right|~\leq~\C~
2^{j}\prod_{i\in\I\cup\J^\flat} 2^{j-\ell_i}.
 \eeq
 Note that $\ell_i\leq j/2+3$ for $i\in\I$ and $\ell_i>j/2+3$ for $i\in\J=\J^\flat\cup\J^\sharp$.
 From (\ref{intersection of cone norm})-(\ref{Lambda norm><}),  we find
\bel{intersection of cone norm*}
\begin{array}{lr}\ds
\left|~\Gamma^\nu_j\cap\Lambda_{\ell j}\cap\left\{2^{j-1}\leq|\xi|<2^{j+1}\right\}~\right|~\leq~\C~
2^{j}2^{|\I|j/2}\prod_{i\in\J^\flat} 2^{j-\ell_i}
\\\\ \ds~~~~~~~~~~~~~~~~~~~~~~~~~~~~~~~~~~~~~~~~~~~~~~~~~~~~~
~=~\C~2^{j}2^{\left(n-1-|\J|\right)j/2}\prod_{i\in\J^\flat} 2^{j-\ell_i}.
\end{array}
\eeq
Recall $\Lambda_\ell$ defined in (\ref{Cone}).
Suppose $\ell_i>j/2+3$ for some $i=1,2,\ldots,n-1$. We have  $|\lambda_i|\leq2^{-\ell_i+1}<2^{-j/2-1}$. Moreover, By definition of $\Gamma_j^\nu$ in (\ref{Gamma}), we have the following observation.
\begin{remark}
Let $\ell_i>{j/2}+3$ for some $i=1,2,\ldots,n-1$. We have
\bel{Inclusion 1}
\Lambda_{\ell }~\subset~\Cup_{\nu~\colon~\xi^\nu_j\in\mathbb{S}^{n-2}_i }~ \Gamma^\nu_j.
\eeq
In particular, if  $\ell_i>{j/2}+3$ for every $i=1,2,\ldots,n-1$, we have
\bel{Inclusion 2}
\Lambda_{\ell }~\subset~\Cup_{\xi^\nu_j=(\tau,\lambda)=(\pm1,0) }~ \Gamma^\nu_j.
\eeq
 \end{remark}
 \begin{figure}[h]
\centering
\includegraphics[scale=0.28]{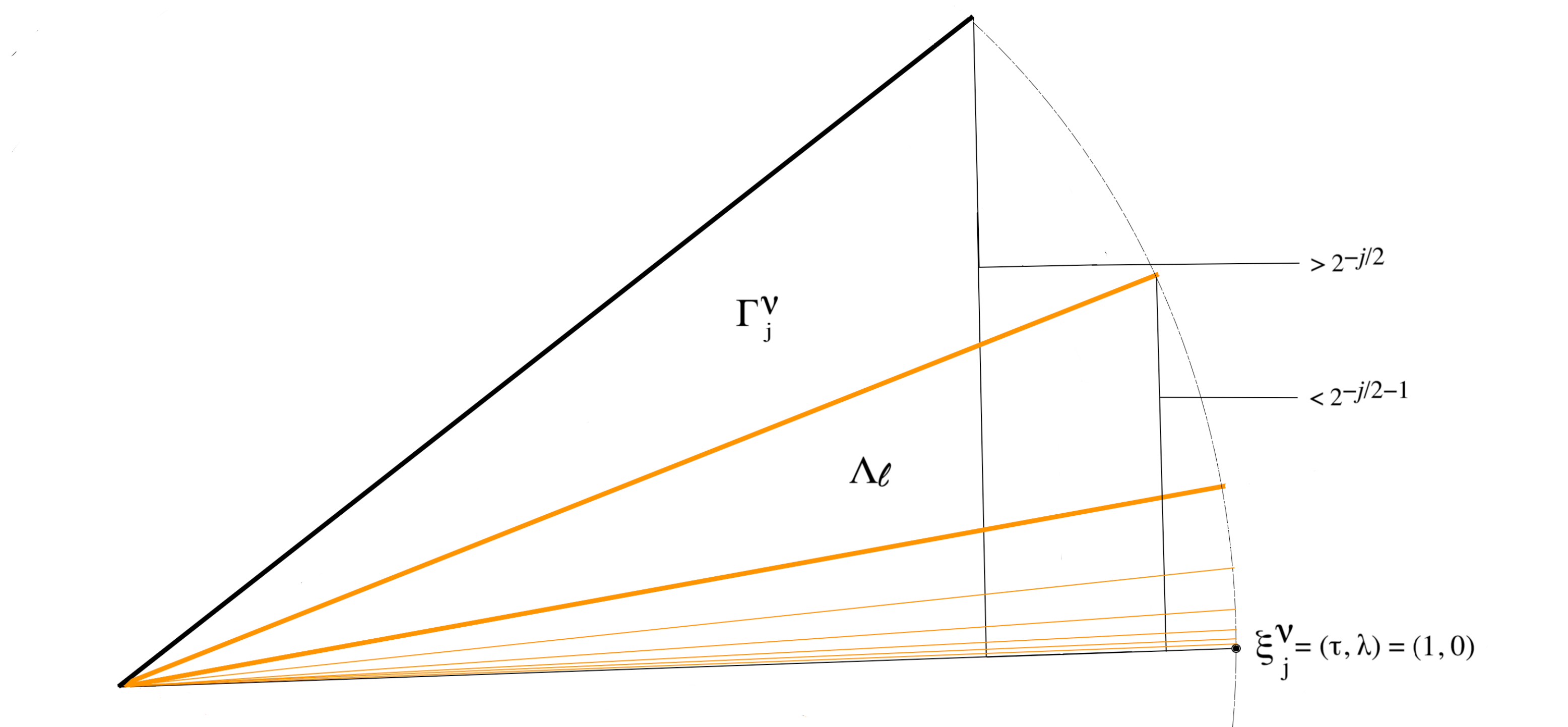}
\caption{\small{$(\tau,\lambda)\in\R\times\R$ and $\ell>j/2+3$.}}
\end{figure}
 
From {\bf Remark 5.3}, if $\J$ is non-empty, $\Lambda_{\ell j}$ can be covered by an union of $\Gamma_j^\nu$ whose central directions belong to $\Cap_{i\in\J}\mathbb{S}^{n-2}_i$.
We define the subset 
\bel{V subset}
\mathcal{V}_{\ell j}~=~\left\{ \nu~\colon ~\xi^\nu_j\in\mathds{S}^{n-1}\cap\Cap_{i\in\J}\mathbb{S}^{n-2}_i,~\Gamma^\nu_j\cap\Lambda_{\ell j}\neq\emptyset\right\}.
\eeq
Let $\varphi$ defined in (\ref{varphi}).  Observe that
\bel{phi^v_j}
\varphi^\nu_j(\xi)~=~\varphi\Bigg[2^{j/2}\left|{\xi\over |\xi|}-\xi^{\nu}_{j}\right|\Bigg]
\eeq
is supported in the geometric cone $\Gamma^\nu_j$.

For every $\nu\in\mathcal{V}_{\ell j}$, we define
 \bel{chi^v_j}
\begin{array}{cc}\ds
\vartheta^\nu_{\ell j}(\xi)~=~\varphi^\nu_j(\xi)\Bigg/\sum_{\mathcal{V}_{\ell j}}\varphi^\nu_j(\xi).
 \end{array}
\eeq
\begin{remark}
Let $\Lambda_{\ell j}$ defined in (\ref{Lambda_lj}).
From {\bf ( 1 )}-{\bf ( 3 )},  there are at most a constant multiple of
$\ds 2^{\left(n-1-|\J|\right)j/2}\prod_{i\in\I} 2^{-\ell_i}$ 
many elements in $\{\xi^{\nu}_j\}_\nu$ such that $\ds\xi^\nu_j\in\mathds{S}^{n-1}\cap\Cap_{i\in\J}\mathbb{S}^{n-2}_i$ and $\Gamma^\nu_j\cap\Lambda_{\ell j}\neq\emptyset$.
\end{remark}
For every $\nu$ fixed, we consider a linear isometry:
$\xi=\L_\nu \eta$ where $\L_\nu$ is an $n\times n$-matrix with $\det\L_\nu=1$. 
In particular,  the $\imath$-th coordinate of $\eta$ is in the same direction of $\xi^\nu_j$ for some $\imath\in\{1,2,\ldots,n\}$. 

Denote $\eta^\nu_j=\left({\eta_\imath\over|\eta_\imath|},0\right)\in\R\times\R^{n-1}$.  We have 
\bel{eta^v_j}
\xi^\nu_j~=~\L_\nu \eta^\nu_j.
\eeq
Furthermore,   we require $\xi_i=\eta_i$ for every $i\in\J$ as $\ds\xi^\nu_j\in\Cap_{j\in\J} \mathds{S}_i^{n-2}$. In the special case of $\J=\{1,2,\ldots,n-1\}$, $\L_\nu$ is the identity matrix so that
$\eta_\imath=\tau$.

\begin{remark}
In the new coordinate system of $\eta\in\R^n$, there is NO definition for $\I=\{1,2,\ldots n-1\}\setminus\J$ such that $0\leq\ell_i\leq j/2+3, i\in\I$.
\end{remark}
Let  $\vartheta^\nu_{\ell j}(\xi)$ defined  in (\ref{chi^v_j}). From direct computation, we find
\bel{chi prop est}
\begin{array}{cc}\ds
\left|\p_\eta^\alpha \vartheta_{\ell j}^{\nu}\left(\L_\nu\eta\right)\right|~\leq~\C_{\alpha} ~2^{|\alphaup| \left({1\over 2}\right)j}|\eta|^{-|\alphaup|}
\end{array}
\eeq
for every multi-index $\alpha$. 

Denote $r=|\xi|=|\eta|$.   For every $\L_\nu\eta=\xi\in\Gamma^\nu_j$, the angle between $\eta$ and $\eta_\imath$ is bounded by $ \arcsin(2\cdot2^{-j/2})$.
By using polar coordinates, we have
\bel{radial derivative}
{\p\over \p\eta_\imath}~=~\left({\p r\over \p \eta_\imath}\right){\p\over \p r}+\O\left(2^{-j/2}\right)\cdot\nabla_{\eta_\imath^\dagger}.
\eeq
Note that 
$\p_r\vartheta^\nu_{\ell j}\equiv0$
 because $\vartheta^\nu_{\ell j}(\xi)=\vartheta^\nu_{\ell j}\left(\L_\nu \eta\right)$  is homogeneous of degree zero in $\eta$. Together with 
(\ref{chi prop est}) and (\ref{radial derivative}), we have
\bel{chi prop est split}
\begin{array}{cc}\ds
\left|\p_{\eta_\imath}^\alpha \vartheta_{\ell j}^{\nu}\left(\L_\nu\eta\right)\right|~\leq~\C_{\alpha} ~|\eta|^{-|\alpha|},
\qquad
\left|\p_{\eta_\imath^\dagger}^\beta \vartheta_{\ell j}^{\nu}\left(\L_\nu\eta\right)\right|~\leq~\C_{\beta} ~2^{|\beta|j/2}|\eta|^{-|\beta|}
\end{array}
\eeq
for every multi-indices $\alpha, \beta$.

\section{Region of influence}
\setcounter{equation}{0}
Recall $\I\cup\J=\{1,2,\ldots,n-1\}$ is defined in (\ref{IJ}).
From the previous section, we have $\eta_j^\nu=\left({\eta_\imath\over |\eta_\imath|},0\right)\in\R\times\R^{n-1}$ for some $\imath\in\{1,2,\ldots,n\}$ such that 
\bel{xi eta containing} 
\xi^\nu_j~=~\L_\nu\eta^\nu_j~\in~\mathds{S}^{n-1}\cap\Cap_{i\in\J} \mathds{S}_i^{n-2}.
\eeq
On the other hand,  $\xi_i=\eta_i$ for every $i\in\J$. Therefore, we must have $\imath\notin\J$.

Consider the rectangle
\bel{rectangle R}
\begin{array}{lr}\ds
R^\nu_j(x_o)~=~
\Bigg\{x\in\supp\sigma~\colon~ 
\left|\left(\L_\nu^Tx_o-\nabla_\eta\Phi\left(x,\L_\nu\eta_j^\nu\right)\right)_\imath\right|\leq4\cdot2^{-j},
\\\\ \ds~~~~~~~~~~~~~~~~
\left\{\sum_{i\neq\imath, i\notin\J}\left(\L_\nu^Tx_o-\nabla_\eta\Phi\left(x,\L_\nu\eta_j^\nu\right)\right)_i^2\right\}^{1\over 2}\leq 4\cdot2^{-j/2}
\Bigg\}.
\end{array}
\eeq
\begin{remark}. There is no restriction for 
\bel{J no restrict}
\left(\L_\nu^Tx_o-\nabla_\eta\Phi\left(x,\L_\nu\eta_j^\nu\right)\right)_i, \qquad  i\in\J.
\eeq
\end{remark}
The set $\mathfrak{Q}_r(x_o)$ is defined by
\bel{Q_r}
\mathfrak{Q}_r(x_o)~\doteq~\bigcup_{2^{-j}\leq r}~\Bigg(~\bigcup_{\nu~\colon~\xi^\nu_j\in\Cap_{i\in\J}\mathbb{S}^{n-2}_i} R_{j}^{\nu}(x_o)~\Bigg).
\eeq
Note that  there are at most a constant multiple of $2^{j\left(n-1-|\J|\right)/2}$ elements in $\{\xi^\nu_j\}_\nu$ for which $\xi^\nu_j\in\mathds{S}^{n-1}\cap\Cap_{i\in\J} \mathbb{S}^{n-2}_i$. 
We have
\bel{Q_r Est exact}
\begin{array}{lr}\ds
 \left|\Q_r(x_o)\right|~\leq~\sum_{2^{-j}\leq r}~~\sum_{\nu~\colon~\xi^\nu_j\in\mathds{S}^{n-1}\cap\Cap_{i\in\J}\mathbb{S}^{n-2}_i}\left|R^{\nu}_{j}(x_o)\right|
 \\\\ \ds ~~~~~~~~~~~~~
 ~\leq~\C_{\sigma~\Phi}~\sum_{2^{-j}\leq r}~~\sum_{\nu~\colon~\xi^\nu_j\in\mathds{S}^{n-1}\cap\Cap_{i\in\J}\mathbb{S}^{n-2}_i} 2^{-j\left(n-1-|\J|\right)/2}2^{-j}\qquad\hbox{\small{by (\ref{rectangle R}) and {\bf Remark 6.2}}}
  \\\\ \ds ~~~~~~~~~~~~~
 ~\leq~\C_{\sigma~\Phi}~\sum_{2^{-j}\leq r} 2^{-j}

  ~\leq~\C_{\sigma~\Phi}~r.
  \end{array}
\eeq

For the associated adjoint operator $\F^*$, we  define 
\bel{rectangle R*}
\begin{array}{lr}\ds
{^*}R^\nu_j(x_o)~=~
\Bigg\{x\in\supp\sigma~\colon~ 
\left|\left(\L_\nu^Tx-\nabla_\eta\Phi\left(x_o,\L_\nu\eta_j^\nu\right)\right)_\imath\right|\leq4\cdot2^{-j},
\\\\ \ds~~~~~~~~~~~~~~~~
\left\{\sum_{i\neq\imath, i\notin\J}\left(\L_\nu^Tx-\nabla_\eta\Phi\left(x_o,\L_\nu\eta_j^\nu\right)\right)_i^2\right\}^{1\over 2}\leq 4\cdot2^{-j/2}
\Bigg\}.
\end{array}
\eeq
whereas $x$ and $x_o$ are switched in  (\ref{rectangle R}).  
The corresponding  region of influence is 
\bel{*Q_r}
\mathfrak{Q}^*_r(x_o)~=~\bigcup_{2^{-j}\leq r}~\Bigg(~\bigcup_{\nu~\colon~\xi^\nu_j\in\mathds{S}^{n-1}\cap\Cap_{i\in\J}\mathbb{S}^{n-2}_i} {^*}R_{j}^{\nu}(x_o)~\Bigg).
\eeq
Clearly,  $\Q^*_r(x_o)$ also satisfies the estimate in (\ref{Q_r Est exact}).
\begin{remark} With all preliminary estimates developed in Section 4 and 5, we are ready to prove {\bf Lemma Two} in the following section. The same argument also applies to $\Omega^*_{\ell j}(x,y)$ defined in (\ref{Omega^*_t,j}) except that $\Q_r(x_o)$ is replaced by $\Q^*_r(x_o)$.
\end{remark}

\section{Proof of Lemma Two}
\setcounter{equation}{0}
Let $\I\cup\J=\{1,2,\ldots,n-1\}$ and $\J=\J^\flat\cup\J^\sharp$  defined in (\ref{IJ}) where
$0\leq\ell_i\leq j/2+3, i\in\I$,  $j/2+3<\ell_i\leq j+3,  i\in\J^\flat$ and $\ell_i>j+3, i\in\J^\sharp$. 

Note that $\xi=\L_\nu\eta$ with $\det\L_\nu=1$. In particular, we have $\xi_i=\eta_i$ for every $i\in\J$. Moreover, there is an $\imath\notin\J$ such that $\ds\eta^\nu_j=\Big(\eta_\imath/|\eta_\imath|,0\Big)\in\R\times\R^{n-1}$ and $\xi^\nu_j=\L_\nu\eta^\nu_j$ as shown in (\ref{eta^v_j}).   

Let $\mathcal{V}_{\ell j}$ defined in (\ref{V subset}) and $\vartheta^\nu_{\ell j}(\xi)$ defined  in (\ref{phi^v_j})-(\ref{chi^v_j}).

Now, recall $\Omega_{\ell j}(x,y)$ from (\ref{Omega_t,j}). 
  We have
\bel{Omega^v_lj}
\begin{array}{cc}\ds
\Omega_{\ell j}(x,y)~=~\sum_{\nu\in\mathcal{V}_{\ell j}} \Omega_{\ell j}^\nu(x,y),
\\\\ \ds
\Omega^\nu_{\ell j}(x,y)~=~\int_{\R^n} e^{2\pi\i\left(\Phi(x,\xi)-y\cdot\xi\right)}\vartheta^\nu_{\ell j}(\xi)\delta_{\ell j}(\xi)\phi_j(\xi)\sigma(x,\xi)d\xi
\end{array}
\eeq
where $\delta_{\ell j}(\xi)$ and $\phi_j(\xi)$ are defined in (\ref{delta_lj}) and (\ref{phi_j}) respectively.

Consider
\bel{Phi split}
\begin{array}{cc}\ds
\Phi(x,\L_\nu\eta)-y\cdot \L_\nu\eta
~=~\left(\nabla_\eta\Phi\left(x,\L_\nu\eta_j^\nu\right)-\L_\nu^T y\right)\cdot\eta~+~\Psi(x,\eta),
\\\\ \ds
\Psi(x,\eta)~\doteq~\Phi(x,\L_\nu\eta)-\nabla_\eta\Phi\left(x,\L_\nu\eta_j^\nu\right)\cdot\eta.
\end{array}
\eeq
We borrow the next result
from {\bf 4.5}, chapter IX of Stein \cite{Stein}:
\bel{d Est Psi}
\left|\p_{\eta_\imath}^\alpha\Psi(x,\eta)\right|~\leq~\C_\alpha~2^{-|\alpha| j},\qquad
\left|\p_{\eta_\imath^\dagger}^\beta\Psi(x,\eta)\right|~\leq~\C_\beta~2^{-|\beta|j/2}
\eeq
for every multi-indices $\alpha, \beta$ whenever $2^{j-1}\leq|\eta|\leq2^{j+1}$.

Rewrite 
\bel{Omega rewrite}
\Omega^\nu_{\ell j}(x,y)~=~\int_{\R^n} e^{2\pi\i\Big(\nabla_\eta\Phi\Big(x,\L_\nu\eta_j^\nu\Big)-\L_\nu^T y\Big)\cdot\eta}\Theta^\nu_{\ell j}(x,\eta)d\eta
\eeq
and
\bel{Theta}
\Theta^\nu_{\ell j}(x,\eta)~=~e^{2\pi\i\Psi(x,\eta)}\vartheta^\nu_{\ell j}(\L_\nu\eta)\delta_{\ell j}(\L_\nu\eta)\phi_j(\L_\nu\eta)\sigma(x,\L_\nu\eta).
\eeq
Observe that 
\bel{Intersection nonempty}
\L_\nu \eta~=~\xi~\in~\Gamma^\nu_j~\cap~\Lambda_\ell~\cap~\left\{2^{j-1}\leq|\xi|=|\eta|\leq2^{j+1}\right\}
\eeq
for $\eta$ in the support of $\Theta^\nu_{\ell j}(x,\eta)$ where $\Gamma^\nu_j$ and $\Lambda_\ell(\xi)$  are defined in (\ref{Gamma}) and  (\ref{Cone}) respectively.

For $\L_\nu \eta\in\Gamma^\nu_j\cap\left\{2^{j-1}\leq|\eta|\leq2^{j+1}\right\}$, we have
\bel{eta norm}
2^{j-1}~\leq~|\eta_\imath|~\leq~2^{j+1},\qquad |\eta_\imath^\dagger|~\leq~\C~ 2^{j/2}.
\eeq
On the other hand, for $\xi=(\tau,\lambda)\in \Lambda_\ell \cap\left\{2^{j-1}\leq|\xi|\leq2^{j+1}\right\}$, we have 
\bel{tau lambda norm}
\begin{array}{cc}\ds
2^{j-1}~\leq~|\tau|~\leq~2^{j+1},
\qquad
 2^{j-1-\ell_i}~\leq~|\lambda_i|~\leq ~2^{j+1-\ell_i},\qquad i=1,2,\ldots,n-1. 
\end{array}
\eeq
Write  $(\tau,\lambda)=\xi=\L_\nu \eta$ for which
\bel{tau lambda entry}
\begin{array}{cc}\ds
\tau~=~a_{ n \imath} \eta_\imath~+~\O(1)\cdot \eta_\imath^\dagger,\qquad ( \xi_n=\tau )
\\\\ \ds
\lambda_i~=~a_{i\imath}\eta_\imath~+~\O(1)\cdot\eta_\imath^\dagger, \qquad i\in\I,\qquad \lambda_i=\eta_i,\qquad i\in\J
\end{array}
\eeq
where  $a_{i\imath}$ denotes the entry on the  $i$-th row and the $\imath$-th column of $\L_\nu$.  

By putting together (\ref{eta norm})-(\ref{tau lambda norm}) and (\ref{tau lambda entry}), we necessarily  have
\bel{entry norm}
\begin{array}{cc}\ds
|a_{n\imath}|~\leq~\C\qquad\hbox{and}\qquad |a_{i\imath}|\leq\C2^{-\ell_i}, \qquad i\in\I.
\end{array}
\eeq
Let $\delta_{\ell j}(\xi)$ defined in (\ref{delta_lj}). Moreover, recall  $\delta_\ell(\xi)$ from (\ref{delta_t})-(\ref{Cone}) satisfying the differential inequality in  (\ref{delta Diff Ineq}). 
Suppose $|\tau|=\C2^j$ and $|\lambda|=\C2^{j-\ell_i}$ as in (\ref{tau lambda norm}).
From  direct computation, for every multi-indices $\alpha,\beta$, we have
\bel{delta_lj Diff Ineq I}
\begin{array}{lr}\ds
\left| \p_\tau^\alpha \prod_{i\in\I}\p_{\lambda_i}^{\beta_i}  \delta_{\ell j}(\tau,\lambda)\right|~\leq~\C_{\alpha~\beta}~ \left({1\over |\tau|}\right)^{\alpha}\prod_{i\in\I} \left({1\over |\lambda_i|}\right)^{\beta_i}
\\\\ \ds~~~~~~~~~~~~~~~~~~~~~~~~~~~~~~~
~\leq~\C_{\alpha~\beta}~2^{-j\alpha} \prod_{i\in\I} 2^{-(j-\ell_i)\beta_i}~\leq~\C_{\alpha~\beta}~2^{-j\alpha} \prod_{i\in\I} 2^{-j\beta_i/2}
\end{array}
\eeq
where $\ell_i\leq j/2+3$ for $i\in\I$, and
\bel{delta_lj Diff Ineq J}
\begin{array}{lr}\ds
\left| \p_\tau^\alpha \prod_{i\in\J^\flat}\prod_{i\in\J^\sharp}\p_{\lambda_i}^{\beta_i}  \delta_{\ell j}(\tau,\lambda)\right|~\leq~\C_{\alpha~\beta}~ \left({1\over |\tau|}\right)^{\alpha}\prod_{i\in\J^\flat} \left({1\over |\lambda_i|}\right)^{\beta_i} \prod_{i\in\J^\sharp} 2^{j\beta_i}\left({1\over |\tau|}\right)^{\beta_i}
\\\\ \ds~~~~~~~~~~~~~~~~~~~~~~~~~~~~~~~~~~~~~
~\leq~\C_{\alpha~\beta}~2^{-j\alpha} \prod_{i\in\J^\flat} 2^{-(j-\ell_i)\beta_i}.
\end{array}
\eeq
Recall $\sigma\in\S^{-{n-1\over 2}}$  satisfying the differential inequality in (\ref{Class}).
Together with (\ref{entry norm})-(\ref{delta_lj Diff Ineq I}),  by using the chain rule of differentiation, we find
\bel{Diff Ineq sigma eta_imath}
\begin{array}{lr}\ds
\left|\p_{\eta_\imath}^N \delta_\ell(\L_\nu\eta)\sigma(x,\L_\nu\eta)\right|~\leq~\C_N~\left({1\over 1+|\eta|}\right)^{n-1\over 2}|\eta|^{-N}
\\\\ \ds~~~~~~~~~~~~~~~~~~~~~~~~~~~~~~~~~~
~\leq~\C_N~ 2^{-j\left({n-1\over 2}\right)}2^{-jN},\qquad N\ge0.
\end{array}
\eeq
Consider $\eta_i$ for $i\notin\J$. From (\ref{Class}) and (\ref{entry norm})-(\ref{delta_lj Diff Ineq I}), by using the chain rule of differentiation,  we find
\bel{Diff Ineq sigma eta_i}
\begin{array}{lr}\ds
\left|\p_ {\eta_i}^{N} \delta_\ell(\L_\nu\eta)\sigma(x,\L_\nu\eta)\right|~\leq~\C_N~
\left({1\over 1+|\eta|}\right)^{n-1\over 2}2^{N\left({1\over 2}\right)j} |\eta|^{-N}
\\\\ \ds~~~~~~~~~~~~~~~~~~~~~~~~~~~~~~~~~~
~\leq~\C_N~ 2^{-j\left({n-1\over 2}\right)}2^{-jN/2}
,\qquad N\ge0.
\end{array}
\eeq
Now, we define the differential operator
\bel{D operator}
\begin{array}{lr}\ds
\mathcal{D}~=~I+2^{2j}\left(\p_{\eta_\imath}\right)^2+2^j\sum_{i\neq\imath,  i\notin\J}\left(\p_{\eta_i}\right)^2+\sum_{i\in\J^\flat} 2^{2(j-\ell_i)}\left(\p_{\eta_i}\right)^2+\sum_{i\in\J^\sharp} \left(\p_{\eta_i}\right)^2
\\\\ \ds~~~
~=~I+2^{2j}\left(\p_{\eta_\imath}\right)^2+2^j\sum_{i\neq\imath,  i\notin\J}\left(\p_{\eta_i}\right)^2+\sum_{i\in\J^\flat} 2^{2(j-\ell_i)}\left(\p_{\xi_i}\right)^2+\sum_{i\in\J^\sharp} \left(\p_{\xi_i}\right)^2.
\end{array}
\eeq
Let $\Theta^\nu_{\ell j}(x,\eta)$  defined in (\ref{Theta}).
From (\ref{Class}), (\ref{chi prop est split}),   (\ref{d Est Psi}),  (\ref{delta_lj Diff Ineq J}) and (\ref{Diff Ineq sigma eta_imath})-(\ref{Diff Ineq sigma eta_i}),  we have
\bel{Theta diff est > sharp}
\left|\mathcal{D}^N \Theta^\nu_{\ell j}(x,\eta)\right|~\leq~\C_N~2^{-j\left({n-1\over 2}\right)},\qquad N\ge0.
\eeq
On the other hand, by using (\ref{intersection of cone norm*}), we have
\bel{Theta support}
\begin{array}{lr}\ds
\left|\supp \Theta^\nu_{\ell j}(x,\eta)\right|~\leq~\C~2^{j}2^{\left(n-1-|\J|\right)j/2}\prod_{i\in\J^\flat} 2^{j-\ell_i}.
\end{array}
\eeq
Recall $\Omega^\nu_{\ell j}(x,y)$ defined in (\ref{Omega rewrite}). Note that $\Omega^\nu_{\ell j}(x,y)$ has a same $x$-compact support of $\sigma(x,\xi)$. From (\ref{Theta diff est > sharp})-(\ref{Theta support}), an $N$-fold integration by parts associated to $\mathcal{D}$ shows that
\bel{Omega rewrite norm > sharp}
\begin{array}{lr}\ds
\left|\Omega^\nu_{\ell j}(x,y)\right|~\leq~\C_{N}~2^{-j\left({n-1\over 2}\right)}~2^{j}2^{\left(n-1-|\J|\right)j/2}\prod_{i\in\J^\flat} 2^{j-\ell_i}
\\\\ \ds
\Bigg\{1+4\pi^2 2^{2j}\left(\nabla_\eta\Phi\left(x,\L_\nu\eta_j^\nu\right)-\L_\nu^T y\right)_\imath^2+4\pi^2 2^{j}\sum_{i\neq\imath, i\notin\J}\left(\nabla_\eta\Phi\left(x,\L_\nu\eta_j^\nu\right)-\L_\nu^T y\right)_i^2 
\\\\ \ds
+4\pi^2 \sum_{i\in\J^\flat}2^{2(j-\ell_i)}\left(\nabla_\eta\Phi\left(x,\L_\nu\eta_j^\nu\right)-\L_\nu^T y\right)_i^2+4\pi^2 \sum_{i\in\J^\sharp}\left(\nabla_\eta\Phi\left(x,\L_\nu\eta_j^\nu\right)-\L_\nu^T y\right)_i^2\Bigg\}^{-N}.
\end{array}
\eeq
Consider a local diffeomorphism  
\bel{local diffeo}
\mathcal{X}_\Phi~\colon~ x~\mt~\left(\L_\nu^T\right)^{-1}\nabla_\eta\Phi\left(x,\L_\nu\eta_j^\nu\right)
\eeq
 whose Jacobian is non-zero provided that $\Phi$ satisfies the non-degeneracy condition  (\ref{nondegeneracy}). 

Denote 
$\mathcal{X}=\mathcal{X}(x)\doteq\left(\L_\nu^T\right)^{-1}\nabla_\eta\Phi\left(x,\L_\nu\eta_j^\nu\right)$.
There are exactly $n-1-|\J|$ many terms in the summation $\ds\sum_{i\neq \imath, i\notin\J}$.
By  using (\ref{Omega rewrite norm > sharp}), we have
\bel{int Omega rewrite norm > sharp}
\begin{array}{lr}\ds
\int_{\R^n}\left|\Omega^\nu_{\ell j}(x,y)\right|dx
~\leq~\C_{\Phi~N}\int_{\R^n} 2^{-j\left({n-1\over 2}\right)}~2^{j}2^{\left(n-1-|\J|\right)j/2}\prod_{i\in\J^\flat} 2^{j-\ell_i}
\\\\ \ds~~~~~~~
\left\{1+ 2^{2j}(\mathcal{X}-y)_\imath^2+2^j\sum_{i\neq\imath, i\notin\J}(\mathcal{X}-y)_i^2+\sum_{i\in\J^\flat}2^{2(j-\ell_i)} (\mathcal{X}-y)_i^2+\sum_{i\in\J^\sharp}(\mathcal{X}-y)_i^2\right\}^{-N}d\mathcal{X} 
\\\\ \ds
~\leq~\C_{\sigma~\Phi~N}\iiiint_{\R\times\R^{n-1-|\J|}\times\R^{|\J^\flat|}\times\R^{|\J^\sharp|}} 2^{-j\left({n-1\over 2}\right)}\left\{1+ \mathcal{Z}_\imath^2+     \sum_{i\neq\imath,i\notin\J} \mathcal{Z}_i^2 +\sum_{i\in\J^\flat}\mathcal{Z}_i^2+\sum_{i\in\J^\sharp} \mathcal{Z}_i^2\right\}^{-N} 
\\\\ \ds~~~~~~~~~~~~~~~~~~~~~~~~~~~~~~~~~~~~~~~~~~~~~~~~~~~~~~~~~~~~~~~~~~~~~~~~~~~~~

d\mathcal{Z}_\imath \prod_{i\neq\imath, i\notin\J}d\mathcal{Z}_i \prod_{i\in\J^\flat}d\mathcal{Z}_i \prod_{i\in\J^\sharp}d\mathcal{Z}_i 
\\\\ \ds~~~~~~~
\hbox{\small{$ \mathcal{Z}_\imath=2^{j}(\mathcal{X}-y)_\imath$,~~ $\mathcal{Z}_i=2^{j/2}(\mathcal{X}-y)_i,~~i\neq\imath, i\notin\J$,~ ~$\mathcal{Z}_i=2^{j-\ell_i}(\mathcal{X}-y)_i,~~ i\in\J^\flat$~and~$\mathcal{Z}_i=(\mathcal{X}-y)_i,~~i\in\J^\sharp$.}}
\\\\ \ds
~\leq~\C_{\sigma~\Phi}~2^{-j\left({n-1\over 2}\right)}\qquad\hbox{\small{for $N$ sufficiently large.}}
\end{array}
\eeq
Recall from {\bf Remark 5.4}.  There are at most $\ds\C ~2^{\left(n-1-|\J|\right)j/2}\prod_{i\in\I} 2^{-\ell_i}$ many elements in $\{\xi^{\nu}_j\}_\nu$ such that $\xi^\nu_j\in\mathds{S}^{n-1}\cap\Cap_{i\in\J}\mathbb{S}^{n-2}_i$ and $\Gamma^\nu_j\cap\Lambda_{\ell j}\neq\emptyset$.
We thus have
\bel{Omega_lj Sum norm > sharp}
\begin{array}{lr}\ds
\int_{\R^n}\left|\Omega_{\ell j}(x,y)\right|dx~\leq~\sum_{\nu~\colon~\xi^\nu_j\in\mathds{S}^{n-1}\cap\Cap_{i\in\J}\mathbb{S}^{n-2}_i, \Gamma^\nu_j\cap\Lambda_{\ell j}\neq\emptyset} \int_{\R^n}\left|\Omega^\nu_{\ell j}(x,y)\right|dx
\qquad \hbox{\small{by (\ref{Omega^v_lj})}}
\\\\ \ds ~~~~~~~~~~~~~~~~~~~~~~~~~~~~
~\leq~\C_{\sigma~\Phi}~ 2^{j\left(n-1-|\J|\right)/2}\prod_{i\in\I}2^{-\ell_i}~2^{-j\left({n-1\over 2}\right)}

\\\\ \ds~~~~~~~~~~~~~~~~~~~~~~~~~~~~
~\leq~\C_{\sigma~\Phi} ~\prod_{i\in\I}2^{-\ell_i} \prod_{i\in\J} 2^{-\left({1\over 2}\right)j}
\\\\ \ds~~~~~~~~~~~~~~~~~~~~~~~~~~~~
~=~\C_{\sigma~\Phi} ~\prod_{i\in\I}2^{-\ell_i} \prod_{i\in\J^\flat} 2^{-\left({1\over 2}\right)j}\prod_{i\in\J^\sharp} 2^{-\left({1\over 2}\right)j}
\\\\ \ds~~~~~~~~~~~~~~~~~~~~~~~~~~~~
~\leq~\C_{\sigma~\Phi} ~\prod_{i\in\I}2^{-\ell_i} \prod_{i\in\J^\flat} 2^{-\left({1\over 2}\right)\ell_i}\prod_{i\in\J^\sharp} 2^{-\left({1\over 2}\right)j}.\qquad (~\ell_i\leq j+3,~i\in\J^\flat~)
\end{array}
\eeq
Observe that every $\p_y$ acting on $\Omega^\nu_{\ell j}(x,y)$ defined in (\ref{Omega^v_lj}) or (\ref{Omega rewrite}) gains a factor of $\C 2^j$ whenever $2^{j-1}\leq|\xi|=|\eta|\leq2^{j+1}$. By carrying out the same estimate in (\ref{Intersection nonempty})-(\ref{Omega_lj Sum norm > sharp}), we find
\bel{nabla Omega_lj>} 
\int_{\R^n} \left|\nabla_y\Omega_{\ell j}(x,y)\right| dx~\leq~\C_{\sigma~\Phi}~
 2^j\prod_{i\in\I}2^{-\ell_i} \prod_{i\in\J^\flat} 2^{-\left({1\over 2}\right)\ell_i}\prod_{i\in\J^\sharp} 2^{-\left({1\over 2}\right)j}.
\eeq
This further implies 
\bel{Est2 >} 
\int_{\R^n} \left|\Omega_{\ell j}(x,y)-\Omega_{\ell j}(x,x_o)\right| dx~\leq~\C_{\sigma~\Phi}~2^j|y-x_o|
 \prod_{i\in\I}2^{-\ell_i} \prod_{i\in\J^\flat} 2^{-\left({1\over 2}\right)\ell_i}\prod_{i\in\J^\sharp} 2^{-\left({1\over 2}\right)j}.
  \eeq

Recall $\mathfrak{Q}_r(x_o)$ defined in (\ref{rectangle R})-(\ref{Q_r}). Let $2^{k}\leq r^{-1}\leq 2^{k+1}$.
For $x\in\R^n\setminus\mathfrak{Q}_r(x_o)$, we either have
\bel{x-Phi Est >}
\left|\left(\L_\nu^Tx_o-\nabla_\eta\Phi\left(x,\L_\nu \eta^\nu_j\right)\right)_\imath\right|~\ge~2\cdot2^{-k}
\eeq
or
\bel{x-Phi Est >'}
\left\{\sum_{i\neq\imath, i\notin\J}\left(\L_\nu^Tx_o-\nabla_\eta\Phi\left(x,\L_\nu\eta_j^\nu\right)\right)_i^2\right\}^{1\over 2}~\ge~ 2\cdot2^{-k/2}.
\eeq
If $y\in B_r(x_o)$, then $|y-x_o|\leq 2^{-k}$. For every $2^j\ge r^{-1}$,  we must have
\bel{x,y-Phi Est >}
\begin{array}{rl}\ds
2^{2j}\left(\L_\nu^Ty-\nabla_\eta\Phi\left(x,\L_\nu \eta^\nu_j\right)\right)_\imath^2+
2^j\sum_{i\neq\imath, i\notin\J} \left(\L_\nu^Ty-\nabla_\eta\Phi\left(x,\L_\nu\eta^\nu_j\right)\right)_i^2
\\\\ \ds
~\ge~2^{2(j-k)}+ 2^{j-k}~\ge~2^{j-k}.
\end{array}
\eeq
Now, repeat the same estimates from (\ref{Intersection nonempty}) to (\ref{Omega_lj Sum norm > sharp}), except that (\ref{Omega rewrite norm > sharp}) is replaced  by the following:

\bel{Omega rewrite norm > sharp*}
\begin{array}{lr}\ds
\left|\Omega^\nu_{\ell j}(x,y)\right|~\leq~\C_{N}~2^{-j\left({n-1\over 2}\right)}~2^{j}2^{\left(n-1-|\J|\right)j/2}\prod_{i\in\J^\flat} 2^{j-\ell_i}
\\\\ \ds
\Bigg\{1+4\pi^2 2^{2j}\left(\nabla_\eta\Phi\left(x,\L_\nu\eta_j^\nu\right)-\L_\nu^T y\right)_\imath^2+4\pi^2 2^{j}\sum_{i\neq\imath, i\notin\J}\left(\nabla_\eta\Phi\left(x,\L_\nu\eta_j^\nu\right)-\L_\nu^T y\right)_i^2 
\\\\ \ds
+4\pi^2 \sum_{i\in\J^\flat}2^{2(j-\ell_i)}\left(\nabla_\eta\Phi\left(x,\L_\nu\eta_j^\nu\right)-\L_\nu^T y\right)_i^2+4\pi^2 \sum_{i\in\J^\sharp}\left(\nabla_\eta\Phi\left(x,\L_\nu\eta_j^\nu\right)-\L_\nu^T y\right)_i^2\Bigg\}^{-N}
\\\\ \ds
~\leq~\C_{N}~2^{-j+k}~2^{-j\left({n-1\over 2}\right)}~2^{j}2^{\left(n-1-|\J|\right)j/2}\prod_{i\in\J^\flat}2^{j-\ell_i}
\\\\ \ds
\Bigg\{1+4\pi^2 2^{2j}\left(\nabla_\eta\Phi\left(x,\L_\nu\eta_j^\nu\right)-\L_\nu^T y\right)_\imath^2+4\pi^2 2^{j}\sum_{i\neq\imath, i\notin\J}\left(\nabla_\eta\Phi\left(x,\L_\nu\eta_j^\nu\right)-\L_\nu^T y\right)_i^2 
\\\\ \ds
+4\pi^2 \sum_{i\in\J^\flat}2^{2(j-\ell_i)}\left(\nabla_\eta\Phi\left(x,\L_\nu\eta_j^\nu\right)-\L_\nu^T y\right)_i^2+4\pi^2 \sum_{i\in\J^\sharp}\left(\nabla_\eta\Phi\left(x,\L_\nu\eta_j^\nu\right)-\L_\nu^T y\right)_i^2\Bigg\}^{1-N}
~~~~
 \hbox{\small{by (\ref{x,y-Phi Est >}).}}
\end{array}
\eeq
 We find
\bel{Est3 >} 
\int_{\R^n\setminus\mathfrak{Q}_r(x_o)} \left|\Omega_{\ell j}(x,y)\right| dx~\leq~\C_{\sigma~\Phi}~ {2^{-j}\over r}~
  \prod_{i\in\I}2^{-\ell_i} \prod_{i\in\J^\flat} 2^{-\left({1\over 2}\right)\ell_i}\prod_{i\in\J^\sharp} 2^{-\left({1\over 2}\right)j} \eeq
for every $y\in B_r(x_o)$ whenever $2^{j}\ge r^{-1}$.

{\small Department of Mathematics, Westlake University}

 {\small email: wangzipeng@westlake.edu.cn}

\end{document}